\documentclass[reqno]{amsart}
\usepackage{epsfig,amsmath,amsfonts,latexsym}
\usepackage{amsthm}
\usepackage[abs]{overpic}
\usepackage[usenames,dvipsnames]{xcolor}
\usepackage{palatino}
\usepackage[hyperfootnotes=false]{hyperref}
\usepackage{caption}
\usepackage{dsfont}
\usepackage{mathtools}
\usepackage{cite}
\usepackage{tikz-cd}
\usepackage{amssymb}
\usepackage{comment}

\newtheorem{fed}{\textbf{Definition}}[section]
\newtheorem{thm}[fed]{\textbf{Theorem}}
\newtheorem{lemma}[fed]{\textbf{Lemma}}

\newtheorem{prop}[fed]{\textbf{Proposition}}

\usepackage{amssymb,graphicx,epsfig,psfrag,epic,eepic,latexsym,color}
\usepackage{amsmath}
\usepackage{mathrsfs}
\usepackage{dsfont}

\newcommand{\R}{\mathbb{R}}

\title{On GIT quotients of the symplectic group,\\stability and bifurcations of symmetric orbits}
\author{Urs Frauenfelder, Agustin Moreno}

\address[Agustin Moreno]{School of Mathematics \\ Institute for Advanced Study \\  Princeton  \\ USA}

\email{agustin.moreno2191@gmail.com}

\address[Urs Frauenfelder]{Institut f\"ur Mathematik \\ Augsburg Universität \\ Augsburg  \\ Germany}

\email{urs.frauenfelder@math.uni-augsburg.de}

\date{}

\begin{document}

\maketitle

\begin{abstract}
    We provide topological obstructions to the existence of orbit cylinders of symmetric orbits, for mechanical systems preserved by antisymplectic involutions (e.g.\ the restricted three-body problem). Such cylinders induce continuous paths which do not cross the \emph{bifurcation locus} of suitable GIT quotients of the symplectic group, which are branched manifolds whose topology provide the desired obstructions. Namely, the complement of the corresponding loci consist of several connected components which we enumerate and explicitly describe; by construction these cannot be joined by a path induced by an orbit cylinder. Our construction extends the notions from Krein theory (which only applies for elliptic orbits), to allow also for the case of symmetric orbits which are hyperbolic. This gives a general theoretical framework for the study of stability and bifurcations of symmetric orbits, with a view towards practical and numerical implementations within the context of space mission design. We shall explore this in upcoming work.
\end{abstract}

\tableofcontents

\section{Introduction} The study of symmetries in classical mechanics is a well-explored and central topic. There are several important mechanical systems which allow for global symmetries in the form of antisymplectic involutions, i.e.\ anti-preserving an ambient symplectic form, but leaving the Hamiltonian invariant. The restricted three-body problem, concerning the dynamics of a negligible mass under the gravitational attraction of two large masses, is a well-known example of such a system. In this setting, a very natural class of objects of study are closed orbits which are symmetric with respect to such an involution. These orbits are also prone to be found by numerical methods, and therefore entail practical interest. For instance, such is the case in the context of space mission design, e.g.\ when placing satellites in orbit around a celestial body, where the stability properties of the orbit in question, as well as knowledge on potential bifurcations, play an important role. 

The list of orbits that have been found, for the three-body problem alone, is certainly long (see e.g.\ \cite[Chapter 9]{Sze} for numerical work, \cite{Henon} for a quantitative analysis of bifurcations, \cite{Bruno} for symmetric planar orbits, \cite{Kalantonis} and references therein for a very recent numerical investigation on the Hill lunar problem); this poses the necessity of keeping track on how they relate to each other. A natural question is then the following: given two symmetric orbits, does there exist a \emph{(symmetric) orbit cylinder} between them, i.e.\ can they be joined by a $1$-parameter family of symmetric orbits which does \emph{not} undergo bifurcation? Alternatively, if a bifurcation is indeed found to be present (e.g.\ by numerical means), can we catalogue it among a finite list of bifurcation types, or can we predict the existence of orbits after a bifurcation by knowledge of orbits before the bifurcation?

In this article, we will provide topological obstructions to the existence of such orbit cylinders, thus addressing the first question. We will also relate them to the stability properties of the corresponding orbits. These obstructions can be cast in terms of properties of the spectrum of the relevant matrices, which can be efficiently implemented in a computer, and therefore used for practical applications.

\medskip

\textbf{Upcoming work.} The second question will be addressed in a separate article, where we consider the notion of the \emph{SFT-Euler characteristic,} as the Euler characteristic of suitable local Floer homology groups, which stays invariant under bifurcations and can be recast in terms which are also amenable for numerical work. Similarly, we consider the \emph{real Euler characteristic}, as the Euler characteristic of the relevant Lagrangian Floer homology groups (arising when the symmetric orbits is thought of as a Lagrangian chord). In practice, these invariants can be used by engineers as a test to predict the existence of orbits: if this number is found to differ before and after the bifurcation, one knows that the algorithm has not found all the orbits. Moreover when combined with the obstructions provided here, they provide educated guesses as to where to look for such orbits, as will explained in upcoming work.

\smallskip

\textbf{Symmetric orbits and monodromy matrices.} A symmetric closed orbit can be thought of as a chord or open string, i.e.\ an orbit segment with its endpoints lying in the fixed-point set of the antisymplectic involution, which is a Lagrangian submanifold of the ambient symplectic manifold. After fixing the energy and projecting out the direction of the flow, the linearization of the dynamics along the orbit gives a time-dependent family of $2n\times 2n$-symplectic matrices (the \emph{reduced monodromy matrices}), all related to each other by symplectic conjugation. If we pick one endpoint of the chord as a starting point, the corresponding matrix at this point satisfies special symmetries. Concretely, such a matrix has the form
\begin{equation}\label{symsymp}
M=M_{A,B,C}=\left(\begin{array}{cc}
A & B\\
C & A^T
\end{array}\right)\in M_{2n\times 2n}(\mathbb{R}),
\end{equation}
where $A,B,C$ are $n\times n$-matrices that satisfy the equations
\begin{equation}\label{eq}
B=B^T,\quad C=C^T,\quad AB=BA^T,\quad
A^TC=CA,\quad A^2-BC=I,
\end{equation}
which ensure that $M$ is symplectic. We will denote the space of such symplectic matrices by $Sp^{\mathcal{I}}(2n)$. 

A natural space to consider is then the quotient space $\mathrm{Sp}(2n)/\mathrm{Sp}(2n)$, where $\mathrm{Sp}(2n)$ acts by conjugation on itself; the above discussion implies that the linearized flow along symmetric orbits admits a specially nice symmetric representative in this quotient. This is a geometric intepretation of the following general algebraic fact due to Wonenburger \cite{wonenburger}: any $2n\times 2n$-symplectic matrix is symplectically conjugate to a symplectic matrix satisfying the above symmetries. 

\smallskip

\textbf{GIT quotients.} In general, there is a slight topological technicality: $\mathrm{Sp}(2n)/\mathrm{Sp}(2n)$ is \emph{not} a Hausdorff space, i.e.\ there are points which cannot be separated from each other. But if one replaces the orbit relation by
the \emph{orbit closure relation}, i.e.\ where two symplectic matrices are identified whenever their orbits under the conjugation action intersect, one obtains the \emph{GIT quotient} $\mathrm{Sp}(2n)//\mathrm{Sp}(2n)$, which does become a Hausdorff space. The transition to the orbit closure relation basically means to ignore Jordan factors, replacing them with diagonal blocks; the resulting matrices, while not necessarily equivalent in $\mathrm{Sp}(2n)/\mathrm{Sp}(2n)$, become so in $\mathrm{Sp}(2n)//\mathrm{Sp}(2n)$. We review this in detail in Appendix \ref{app:GIT} below, and consider only GIT quotients in what follows.  

\smallskip

\textbf{The GIT sequence.} Note that the above expression for $M_{A,B,C}$ implies the choice of a basis for the tangent space to the fixed-point locus along the endpoint of the chord. A different choice of basis amounts to acting with an invertible matrix $R \in GL_n(\mathbb{R})$, via

\begin{equation}\label{act}
R_*\big(A,B,C\big)=\Big(RAR^{-1},RBR^T,(R^T)^{-1}CR^{-1}\Big),
\end{equation}
i.e.\ $M_{A,B,C}$ is replaced by $M_{R_*(A,B,C)}$. Note that this action acts on $A$ by conjugation, and it is also easy to check that $M_{A,B,C}$ and $M_{R_*(A,B,C)}$ are symplectically conjugated. We may then consider the sequence of maps between GIT quotients
\begin{equation}\label{seq}
\mathrm{Sp}^\mathcal{I}(2n)//\mathrm{GL}_n(\mathbb{R}) \mapsto \mathrm{Sp}(2n)//\mathrm{Sp}(2n)
\mapsto \mathrm{M}_{n\times n}(\mathbb{R})//\mathrm{GL}_n(\mathbb{R})
\end{equation}
given by
$$[M_{A,B,C}] \mapsto [[M_{A,B,C}]] \mapsto [A].$$
Here we denote by $[M_{A,B,C}]$ the equivalence class of the matrix
$M_{A,B,C} \in \mathrm{Sp}^\mathcal{I}(2n)$ in the GIT quotient $\mathrm{Sp}^\mathcal{I}(2n)//\mathrm{GL}_n(\mathbb{R})$, by $[[M_{A,B,C}]]$ the equivalence class in the GIT quotient $\mathrm{Sp}(2n)//\mathrm{Sp}(2n)$, and by $[A]$ the equivalence class of the first block $A\in M_{n\times n}(\mathbb{R})$ in $ \mathrm{M}_{n\times n}(\mathbb{R})//\mathrm{GL}_n(\mathbb{R})$, where $\mathrm{GL}_n(\mathbb{R})$ acts on $\mathrm{M}_{n \times n}(\mathbb{R})$ by conjugation. We remark that by mapping the equivalence class of a matrix $A \in \mathrm{M}_{n\times n}
(\mathbb{R})$ to the coefficients of its characteristic polynomial, we get an identification $\mathrm{M}_{n\times n}(\mathbb{R})//\mathrm{GL}_n(\mathbb{R}) \cong \mathbb{R}^n.$ We shall review this nice fact in Appendix \ref{app:GIT}. 

\smallskip

\textbf{Labeled branched manifolds, and normal forms.} In this article, we will explicitly describe the spaces and maps which appear in the above GIT sequence, in the cases $2n=2$ and $2n=4$, although a similar description of this space works in any dimension. Our main focus to the two and four dimensional cases, in which the reduced monodromy matrices are respectively $2\times 2$ and $4\times 4$ symplectic matrices, comes from the fact that these are physically the most meaningful. For instance, such is the case of, respectively, the planar and spatial restricted three-body problem. We will obtain especially nice descriptions of the GIT quotients appearing in the GIT sequence as \emph{labeled} $n$-dimensional branched manifolds (LBMs), in such a way that the maps in the GIT sequence preserve this structure. Namely, (some of) the branches of these LBMs are equipped with positive/negative labels, keeping track of information attached to the eigenvalues of the corresponding matrices. Crucially, this data stays invariant in the presence of an orbit cylinder. Moreover, the branching locus and bifurcation locus in the base $M_{n\times n}(\mathbb{R})//GL_n(\mathbb{R})=\mathbb{R}^n$ of the GIT sequence can be explicitly understood; see Figure \ref{fig:maps} ($n=1$) and Figure \ref{fig:branchlocus} ($n=2$). The resulting diagram for $n=2$, as we learned after rediscovering it in the context of the above GIT sequence, was originally introduced by Broucke in \cite{Broucke}, and it is sometimes referred to as ``Broucke's stability diagram'' in the engineering literature (see also Howard--Mackay \cite{HM} for the cases $n=3,4$, also incorporating the notion of Krein signature for the study of linear stability). In this article, we will refine these studies by incorporating the notion of the \emph{$B$-signature} in the picture, as explained in the following paragraphs. We will also provide explicit normal forms for every equivalence class in these GIT quotients.  

\smallskip

\textbf{Comparison with Krein theory.} These labels are to be understood in the spirit of Krein theory, which roughly speaking is a refinement of the spectrum of a given symplectic matrix, also equipping the eigenvalues with suitable signs. This data completely characterizes the strong stability of the matrix, as proved independently by Krein and Moser; we shall review this in Appendix \ref{app:Krein} (Theorem \ref{Kreinthm}). However, Krein theory only applies for elliptic orbits. Here, we will extend this theory to allow also for matrices of the form $M_{A,B,C}$ with \emph{hyperbolic} eigenvalues, via the notion of \emph{B-positivity/negativity} (Definition \ref{def:Bpos}). This coincides with the notion of Krein-positivity/negativity in the case of elliptic symmetric orbits (Lemma \ref{lemma:BposKrein}). This has the advantage of incorporating Krein theory in a much simpler and efficient way, for the purposes of practical implementations, as is rather straightforward to check whether a matrix is $B$-positive/negative (provided it is presented in its symmetric form $M_{A,B,C}$). Also, this notion is motivated from the theory of symmetric spaces, i.e.\ $Sp^\mathcal{I}(2n)$ can be identified with the space of linear antisymplectic involutions, a symmetric space where the product of two such involutions is given by conjugating one with the other. In this vein, the information carried by the LBM $\mathrm{Sp}(2n)//\mathrm{Sp}(2n)$ coincides precisely with the information carried by Krein theory; its labels only apply for the elliptic case. On the other hand, the information carried by the LBM $\mathrm{Sp}^\mathcal{I}(2n)//\mathrm{GL}_n(\mathbb{R})$ is more refined, as it allows to distinguish more matrices via the associated labels, which apply also to the hyperbolic case. Topologically, this means that $\mathrm{Sp}^\mathcal{I}(2n)//\mathrm{GL}_n(\mathbb{R})$ has more branches than $\mathrm{Sp}(2n)//\mathrm{Sp}(2n)$, some of which get collapsed under the natural map in the GIT sequence; see Figure \ref{fig:bifurcations}.  

\smallskip

\textbf{The topological obstructions.} The topology of the spaces $\mathrm{Sp}^\mathcal{I}(2n)//\mathrm{GL}_n(\mathbb{R})$ and $\mathrm{Sp}(2n)//\mathrm{Sp}(2n)$, together with the labels, provide precisely the obstruction to the existence of orbit cylinders. Namely, removing the \emph{bifurcation locus} from each of them (consisting of matrices with $\pm 1$ as an eigenvalue, and hence corresponding to bifurcation/period doubling) results in two LBMs with several connected components. The maps in the GIT sequence preserve the bifurcation locus, and map connected components to connected components, acting as covering maps away from the branching locus, with varying covering degree. Matrices in different components, by construction, cannot be connected to each other by a continuous path, hence obstructing the existence of orbit cylinders in the case the matrices arise by linearization along symmetric orbits. In fact, the complement of the bifurcation locus in $\mathrm{Sp}(2n)//\mathrm{Sp}(2n)$ consists precisely of $8$ connected components, whereas its complement in $\mathrm{Sp}^\mathcal{I}(2n)//\mathrm{GL}_n(\mathbb{R})$, of $19$ connected components. This illustrates, in a quantitative way, how much more information is carried by $\mathrm{Sp}^\mathcal{I}(2n)//\mathrm{GL}_n(\mathbb{R})$ when compared to $\mathrm{Sp}(2n)//\mathrm{Sp}(2n)$. In what follows, we carry out the details of this construction. 

\medskip

\textbf{Acknowledgements.} The second author acknowledges the support by the National Science Foundation under Grant No. DMS-1926686.

\section{Geometric and dynamical setup}

We now describe the general setup, which motivates the linear algebra of the sections to come. We assume that $(M,\omega)$ is a symplectic manifold and $H \colon M \to \mathbb{R}$ a smooth function. The Hamiltonian vector field $X_H$ of $H$ is
defined by the requirement that
$$dH=\omega(\cdot,X_H).$$
Abbreviate by $S^1=\mathbb{R}/\mathbb{Z}$ the circle. A periodic orbit 
$x \in C^\infty(S^1,M)$ is a solution of the ODE
$$\partial_t x(t)=\tau X_H(x(t)), \quad t \in S^1,$$
where $\tau$ is a positive real number referred to as the period of the periodic orbit. 
If $\phi^t_H$ denotes the flow of the Hamiltonian vector field of $H$, we can characterize
the periodic orbit equivalently by 
$$x(t)=\phi^{\tau t}_H(x(0)).$$
Abbreviating $x_0=x(0)$ the differential 
$$d\phi^\tau(x_0) \colon T_{x_0}M \to T_{x_0}M$$
is a linear symplectic map of the symplectic vector space
$(T_{x_0}M,\omega_{x_0})$ called the \emph{monodromy}. 
\\ \\
We assume now that the periodic orbit $x$ is nonconstant which is equivalent to the requirement
that $X_H(x(t))$ is never zero or in other words $x(t)$ is no critical point of $H$ for every
$t \in S^1$. Since $H$ is autonomous, i.e., does not depend on time, it follows that
$$d\phi^\tau_H(x_0) X_H(x_0)=X_H(x_0),$$
i.e., $X_H(x_0)$ is an eigenvector to the eigenvalue $1$ of the monodromy. Moreover, 
by preservation of energy $H$ is constant along the periodic orbit and therefore
the monodromy maps the codimension one subspace $\ker dH(x_0)$ of the tangent space
$T_{x_0}M$ into itself. Therefore the monodromy induces a linear symplectic map
on the quotient space
$$\overline{d\phi^\tau_H(x_0)} \colon \ker dH(x_0)/\langle X_H(x_0)\rangle \to
\ker dH(x_0)/\langle X_H(x_0)\rangle$$
which we refer to as the \emph{reduced monodromy}. Hence if the dimension of the symplectic
manifold is $2n$ we can associate to the periodic orbit $x$ an element in
$\mathrm{Sp}(2n-2)//\mathrm{Sp}(2n-2)$, namely the equivalence class of its reduced monodromy.
It is interesting to remark that this class does not depend on the starting point of the periodic
orbit. In fact, if we translate our periodic orbit in time
$$r_*x(t)=x(t+r), \quad t\in S^1$$
for $r \in S^1$, we obtain different (parametrized) periodic orbits $r_* x$ whose reduced
monodromy gives rise to the same element in $\mathrm{Sp}(2n-2)//\mathrm{Sp}(2n-2)$, since
the reduced monodromies at different starting points of the periodic orbit are symplectically conjugated to each other via the differential of the flow of the Hamiltonian vector field.  
\\ \\
We now consider a \emph{real symplectic manifold} $(M,\omega,\rho)$. This is a symplectic manifold
$(M,\omega)$ together with an antisymplectic involution $\rho$, i.e., a diffeomorphism of $M$ 
satisfying
$$\rho^2=\mathrm{id}_M, \quad \rho^*\omega=-\omega.$$
We assume that $H$ is invariant under $\rho$, i.e.,
$$H \circ \rho=H.$$
This implies that the Hamiltonian vector field is antiinvariant, meaning that
$$\rho^*X_H=-X_H.$$
In particular, we obtain for its flow
\begin{equation}\label{anti}
\rho \phi^t_H \rho=\phi^{-t}_H.
\end{equation}
A periodic orbit $x$ of $H$ is called \emph{symmetric} if it satisfies
$$x(t)=\rho(x(-t)), \quad t \in S^1.$$
In particular, we have for a symmetric periodic orbit that
$$x\big(0\big),\,\,x\big(\tfrac{1}{2}\big) \in \mathrm{Fix}(\rho).$$
The fixed point set of an antisymplectic involution $\mathrm{Fix}(\rho)$ is a Lagrangian submanifold
of $M$, and in particular, if we look just at half the symmetric periodic orbit we get a chord
$\mathrm{Fix}(\rho)$ to itself. Hence we can think of a symmetric periodic orbit in two ways,
either as a closed string or as an open string from the Lagrangian $\mathrm{Fix}(\rho)$ to itself. 
\\ \\
The differential of the antisymplectic involution
$$d\rho(x_0) \colon T_{x_0}M \to T_{x_0}M$$
gives rise to a linear antisymplectic involution on the symplectic vector space
$(T_{x_0}M,\omega_{x_0})$ which induces a linear antisymplectic involution on the quotient space
$$\overline{d\rho(x_0)} \colon \ker dH(x_0)/\langle X_H(x_0)\rangle \to
\ker dH(x_0)/\langle X_H(x_0)\rangle.$$
Differentiating (\ref{anti}) we get
$$\overline{d \rho(x_0)} \circ \overline{d\phi^\tau_H(x_0)}
\circ \overline{d \rho(x_0)}
=\big(\overline{d\phi^\tau_H(x_0)}\big)^{-1},$$
i.e., equation (\ref{conj}) for
$$\mathcal{I}=\overline{d \rho(x_0)}, \quad M=\overline{d\phi^\tau_H(x_0)}.$$
Hence to a symmetric periodic orbit the reduced monodromy associates an element in
$\mathrm{Sp}^\mathcal{I}(2n-2)//\mathrm{GL}_{2n-2}(\mathbb{R})$. To obtain this map
we have to choose the starting point of the periodic orbit on the Lagrangian $\mathrm{Fix}(\rho)$.

\section{The symplectic group, symmetries, and GIT quotients} Although we later restrict to dimension four we start our discussion for the general case. Our starting point is a fascinating theorem due to Wonenburger \cite{wonenburger} which tells us that every element
$M \in \mathrm{Sp}(2n)$ can be written as the product of two linear antisymplectic involutions
$$M=\mathcal{I}_1 \mathcal{I}_2.$$
Since $\mathcal{I}_1$ and $\mathcal{I}_2$ are involutions, it follows that
$$M^{-1}=\mathcal{I}_2 \mathcal{I}_1=\mathcal{I}_1 M \mathcal{I}_1,$$
i.e. $M$ is conjugated to its inverse via an antisymplectic involution. All linear antisymplectic
involutions are conjugated to each other, and in particular to the standard antisymplectic involution
$$\mathcal{I}=\left(\begin{array}{cc}
I & 0\\
0 & -I
\end{array}\right)$$
where $I$ is the $n\times n$ identity matrix. Hence there exists $G \in \mathrm{Sp}(2n)$ such that
$$G^{-1}\mathcal{I}_1 G=\mathcal{I}$$
and therefore
$$G^{-1}MG=G^{-1} \mathcal{I}_1G G^{-1}\mathcal{I}_2 G=\mathcal{I}\mathcal{I}_2^G,$$
where $\mathcal{I}_2^G=G^{-1}\mathcal{I}_2 G$ is itself an antisymplectic involution. Hence after
conjugation we can assume that
\begin{equation}\label{conj}
M^{-1}=\mathcal{I}M \mathcal{I}.
\end{equation}
If we write $M$ as a block matrix
$$M=\left(\begin{array}{cc}
A & B\\
C & D
\end{array}\right)$$
for $n\times n$-matrices $A,B,C,$ and $D$, it follows since $M$ is symplectic that these
matrices satisfy
$$AB^T=BA^T,\quad CD^T=DC^T, \quad AD^T-BC^T=I.$$
Moreover, the inverse of $M$ is given by
$$M^{-1}=\left(\begin{array}{cc}
D^T & -B^T\\
-C^T& A^T
\end{array}\right).$$
It follows from (\ref{conj}) that
$$\left(\begin{array}{cc}
D^T & -B^T\\
-C^T& A^T
\end{array}\right)=\left(\begin{array}{cc}
A & -B\\
-C & D
\end{array}\right).$$
Therefore $D=A^T$, and so $M$ can be written as
\begin{equation}\label{symsymp}
M=\left(\begin{array}{cc}
A & B\\
C & A^T
\end{array}\right)
\end{equation}
where $A,B,C$ satisfy the equations
\begin{equation}\label{eq}
B=B^T,\quad C=C^T,\quad AB=BA^T,\quad
A^TC=CA,\quad A^2-BC=I.
\end{equation}
We summarize this discussion in the following proposition.
\begin{prop}
Every symplectic matrix $M \in \mathrm{Sp}(2n)$ is symplectically conjugated to a matrix of the form (\ref{symsymp})
with $A,B,C$ satisfying (\ref{eq}). 
\end{prop}
As the above discussion shows, symplectic matrices of the from (\ref{symsymp}) with
$A,B,C$ satisfying (\ref{eq}) are in one-to-one correspondence with symplectic matrices
$M$ satisfying (\ref{conj}). We abbreviate this submanifold of $\mathrm{Sp}(2n)$ by
$$\mathrm{Sp}^\mathcal{I}(2n)=\big\{M \in \mathrm{Sp}(2n): M^{-1}=\mathcal{I}M \mathcal{I}\big\}.$$
This space itself has an interesting structure. Note that it follows from (\ref{conj}) that
$$(\mathcal{I}M)^2=\mathcal{I}M \mathcal{I}M=M^{-1}M=I$$
so that $\mathcal{I}M$ is itself an antisymplectic involution. Therefore
$M=\mathcal{I}(\mathcal{I}M)$ is precisely a Wonenburger decomposition of $M$ into the product
of two antisymplectic involutions. Therefore one can identify the space
$\mathrm{Sp}^{\mathcal{I}}(2n)$ via the map $M \mapsto \mathcal{I}M$ with the space of
linear antisymplectic involutions, which itself corresponds to the tangent bundle of the
Lagrangian Grassmannian \cite{albers-frauenfelder}. 
\\ \\
In the following we will freely identify the space $\mathrm{Sp}^\mathcal{I}(2n)$ as the moduli space
$$\mathrm{Sp}^\mathcal{I}(2n)=\Big\{(A,B,C) \in \mathrm{M}_{n\times n}(\mathbb{R}): (A,B,C)\,\,\textrm{solution of (\ref{eq})}\Big\},$$
via the map
$$(A,B,C) \mapsto M_{A,B,C}:=\left(\begin{array}{cc}
A & B\\
C & A^T
\end{array}\right).$$
Since every symplectic matrix is symplectically conjugated to one in $\mathrm{Sp}^\mathcal{I}(2n)$ it suffices
to restrict one's attention to this submanifold in order to understand the 
GIT quotient $\mathrm{Sp}(2n)//\mathrm{Sp}(2n)$. Although it might be in general cumbersome to find
for a general matrix $M \in \mathrm{Sp}(2n)$ a matrix conjugated to $M$ in $\mathrm{Sp}^\mathcal{I}(2n)$,
we explain that for reduced monodromy matrices of symmetric periodic orbits in the restricted three-body problem there is a simple geometric way to do that. The message of this paper which we want to convey is that
given a symmetric form (\ref{symsymp}) of a symplectic matrix in its similarity class it is advantageous
to keep it for further exploration of the similarity class. A first hint of this philosophy is provided
by the following lemma, which tells us that the characteristic polynomial $p_{A,B,C}$ of the symplectic matrix
$M_{A,B,C}$ is completely determined by the matrix $A$ and does not depend on the matrices $B$ and $C$.
\begin{lemma}\label{cha}
The characteristic polynomial of a matrix $M_{A,B,C} \in \mathrm{Sp}^\mathcal{I}(2n)$ is given by
$$p_{A,B,C}(t)=t^n p_{-2A}\big(-t-\tfrac{1}{t}\big),$$
where $p_{-2A}$ is the characteristic polynomial of the matrix $-2A$.
\end{lemma}
We shall prove this lemma in Section \ref{sec:charpol} below. Now, the group $\mathrm{GL}_n(\mathbb{R})$ acts on $\mathrm{Sp}^\mathcal{I}(2n)$ by
\begin{equation}\label{act}
R_*\big(A,B,C\big)=\Big(RAR^{-1},RBR^T,(R^T)^{-1}CR^{-1}\Big),
\end{equation}
where $R \in \mathrm{GL}_n(\mathbb{R})$ and $(A,B,C) \in \mathrm{Sp}^\mathcal{I}(2n)$. As explained in the Introduction, this action comes from the ambiguity in choosing a basis for the tangent space to the Lagrangian fixed-point locus, at an endpoint of a chord. Note that $A$ transforms as a linear map, whereas $B,C$ transform as bilinear forms. This corresponds to the conjugation action by a linear symplectomorphism:
$$M_{R_*(A,B,C)}=\left(\begin{array}{cc}
R & 0\\
0 & (R^T)^{-1}
\end{array}\right)\left(\begin{array}{cc}
A & B\\
C & A^T
\end{array}\right)\left(\begin{array}{cc}
R^{-1} & 0\\
0 & R^T
\end{array}\right).$$
We therefore obtain a sequence of maps between GIT quotients
\begin{equation}\label{seq}
\mathrm{Sp}^\mathcal{I}(2n)//\mathrm{GL}_n(\mathbb{R}) \mapsto \mathrm{Sp}(2n)//\mathrm{Sp}(2n)
\mapsto \mathrm{M}_{n\times n}(\mathbb{R})//\mathrm{GL}_n(\mathbb{R})
\end{equation}
given by
$$[M_{A,B,C}] \mapsto [[M_{A,B,C}]] \mapsto [A],$$ 
as explained in the Introduction.

\smallskip

\textbf{Warm-up: two dimensional case.} Let us first describe the sequence (\ref{seq}) in the simplest possible case, i.e.\ $n=1$. This case has also been studied in \cite[Appendix B]{BZ}, where it plays an important role when trying to define a real version of ECH. The identification 
$$\mathrm{M}_{1\times 1}(\mathbb{R})//\mathrm{GL}_1(\mathbb{R}) \cong \mathbb{R}$$
is tautological, and the relevant maps are $[M_{A,B,C}]\mapsto [[M_{A,B,C}]]\mapsto [A]=A=\mathrm{tr}(M_{A,B,C})/2$. 
\begin{figure}
    \centering
    \includegraphics[width=0.37\linewidth]{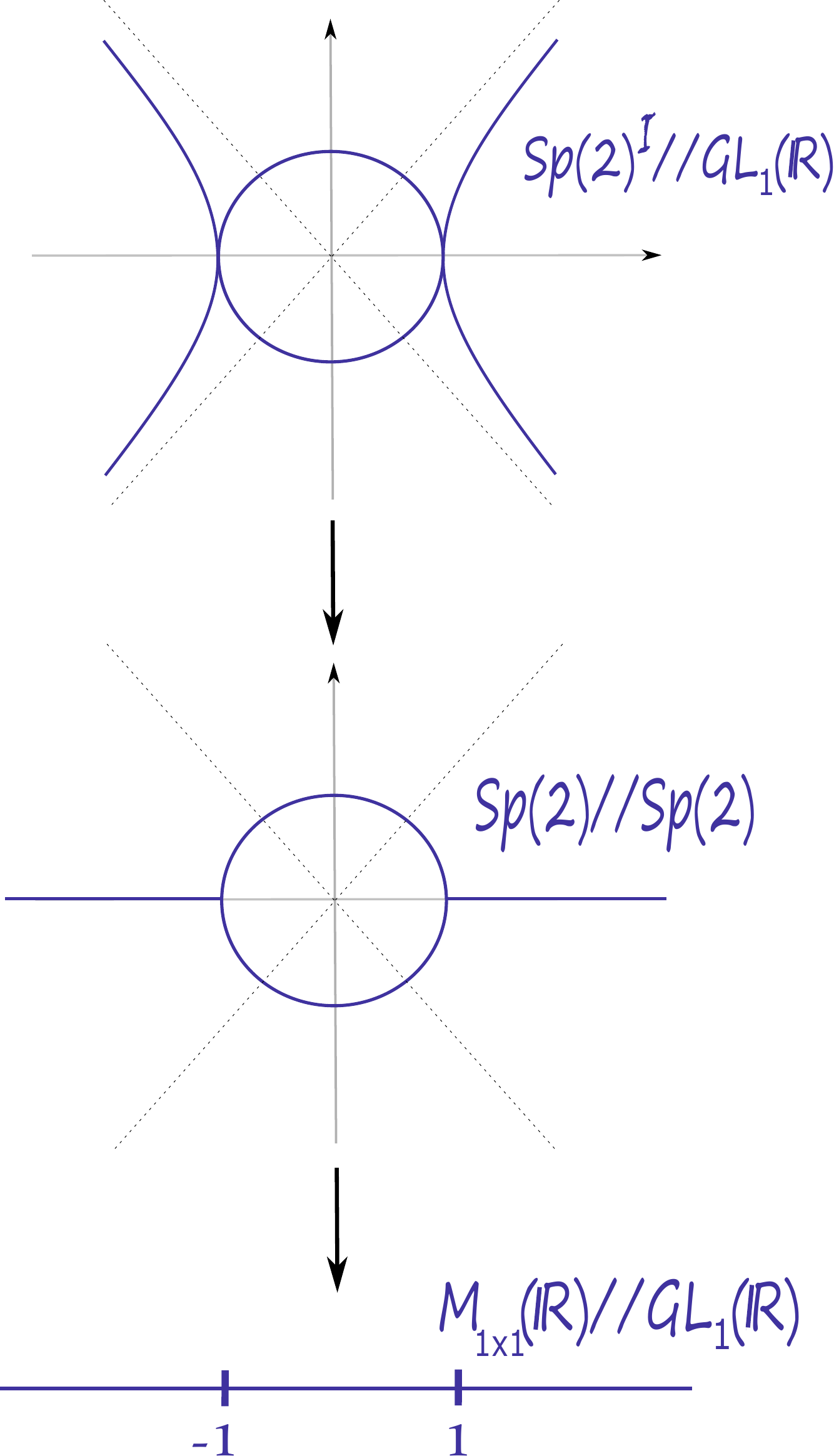}
    \caption{The GIT sequence in the case $n=1$.}
    \label{fig:maps}
\end{figure}
The action of $GL_1(\mathbb{R})=\mathbb{R}^+$ on $\mathrm{Sp}^\mathcal{I}(2)$ is simply 
$$
\epsilon \cdot \left(\begin{array}{cc}
    A & B \\
    C & A
\end{array} \right)
=\left(\begin{array}{cc}
    A & \epsilon^2 B \\
    \frac{1}{\epsilon^2}C & A
\end{array} \right),$$ where $A^2-BC=1$, $\epsilon>0$. We have $\mathrm{Sp}(2)=SL(2,\mathbb{R})$, and a matrix $A\in \mathrm{Sp}(2)$ is either hyperbolic (i.e.\ $\vert \mathrm{tr}(A) \vert >2$, in which case it has two real eigenvalues $r,1/r$ with $\vert r \vert >1$), elliptic (i.e.\ $\vert \mathrm{tr}(A) \vert <2$, in which case it has two conjugate complex eigenvalues in the unit circle), or parabollic (i.e.\ $\vert \mathrm{tr}(A) \vert =2$, in which case it has eigenvalue $\pm 1$ with algebraic multiplicity two). From the discussion in \cite[Section 10.5]{FvK}, we gather that $\mathrm{Sp}(2)//\mathrm{Sp}(2)$ admits a homeomorphism 
$$
\mathrm{Sp}(2)//\mathrm{Sp}(2)=\{z\in \mathbb{C}: \vert z\vert =1\}\cup \{r\in \mathbb{R}: \vert r\vert\geq 1\}\subset \mathbb{C},
$$
via the identification
$$
s(e^{i\theta})=\left( \begin{array}{cc}
\cos(\theta)     & -\sin(\theta) \\
\sin(\theta)     & \cos(\theta)
\end{array}\right), \; s(r)= \left( \begin{array}{cc}
r    & 0 \\
0     & \frac{1}{r}
\end{array}\right).
$$
The hyperbolic locus consists of closed orbits and corresponds to $\{\vert r \vert >1\}$; the elliptic locus also consists of closed orbits, and corresponds to $\{\vert z\vert=1\}\backslash\{\pm 1\}$; and the parabollic locus is $\{\pm 1\}$, where $\{+1\}$ corresponds to the three different Jordan forms with eigenvalue $1$ of algebraic multiplicity two, and, similarly $\{-1\}$ corresponds to the three Jordan forms with eigenvalue $-1$ of algebraic multiplicity two.

Similarly, the GIT quotient $\mathrm{Sp}^\mathcal{I}(2)//\mathrm{GL}_1(\mathbb{R})$ admits an identification
$$
\mathrm{Sp}^\mathcal{I}(2)//\mathrm{GL}_1(\mathbb{R})=\{z\in \mathbb{C}:\vert z\vert=1\}\cup \{(\pm \cosh(u),\sinh(u)):u\in \mathbb{R}\}\subset \mathbb{C},
$$
via
$$
t(e^{i\theta})=s(e^{i\theta})=\left( \begin{array}{cc}
\cos(\theta)     & -\sin(\theta) \\
\sin(\theta)     & \cos(\theta)
\end{array}\right), \; t(u)= \left( \begin{array}{cc}
\pm \cosh(u)    & \sinh(u) \\
\sinh(u)     & \pm \cosh(u)
\end{array}\right).
$$
The matrix $t(u)$ has eigenvalues $\pm e^u, \pm e^{-u}$. Moreover, the matrices $t(u)$ and $t(-u)$ are both symplectically conjugate to diag$(\pm e^u,\pm e^{-u})$, hence to each other, and therefore define the same element in $\mathrm{Sp}(2)//\mathrm{Sp}(2)$.
After these identifications, the GIT sequence becomes
$$
e^{i\theta}\mapsto e^{i\theta} \mapsto \cos(\theta),
$$
$$
(\pm \cosh(u),\sinh(u))\mapsto r=\pm e^{\vert u \vert} \mapsto r=\pm e^{\vert u \vert}.
$$
This sequence is topologically depicted in Figure \ref{fig:maps}; note that it consists of branched maps between $1$-dimensional branched manifolds, with branching locus $\{\pm 1\}$. The covering degree of the first map is two over the hyperbolic locus, and one everywhere else. For the second map, it is two over the elliptic locus, and one elsewhere. 
\smallskip

\textbf{Four dimensional case.} We now describe the sequence (\ref{seq}) in the four dimensional case, where $n=2$. In this
case the identification
$$\mathrm{M}_{2\times 2}(\mathbb{R})//\mathrm{GL}_2(\mathbb{R}) \cong \mathbb{R}^2$$
is obtained via the map
$$\mathrm{M}_{2\times 2}(\mathbb{R})//\mathrm{GL}_2(\mathbb{R}) \to \mathbb{R}^2, \quad
[A] \mapsto \big(\mathrm{tr}(A),\det(A)\big).$$
The spaces $\mathrm{Sp}^\mathcal{I}(4)//\mathrm{GL}_2(\mathbb{R})$ and
$\mathrm{Sp}(4)//\mathrm{Sp}(4)$ are not manifolds but have some branch points. The branch
points lie over three curves in $\mathrm{M}_{2\times 2}(\mathbb{R})//\mathrm{GL}_2(\mathbb{R}) =
 \mathbb{R}^2$ which we describe next. We abbreviate coordinates on $\mathbb{R}^2$ by
 $(\tau,\delta)$, where $\tau$ stands for trace and $\delta$ for determinant. The first branch
 locus is the graph of the parabola $\delta=\tfrac{1}{4}\tau^2$
 $$\Gamma_{d}=\Big\{\big(\tau,\tfrac{1}{4}\tau^2\big): \tau \in \mathbb{R}\Big\}$$
 at which the characteristic polynomial of $A$ has a double root, and the two other branch loci
 are the straight lines
$$\Gamma_1=\Big\{(\tau,\tau-1):\tau \in \mathbb{R}\Big\}, \quad
 \Gamma_{-1}=\Big\{(\tau,-\tau-1):\tau \in \mathbb{R}\Big\}$$  
at which the characteristic polynomial has a root at $1$ respectively $-1$, i.e., the matrix
$A$ has $1$ respectively $-1$ as an eigenvalue. We shall refer to the preimage of $\Gamma_1\cup \Gamma_{-1}$ under the map $\mathrm{Sp}(4)//\mathrm{Sp}(4)\rightarrow \mathbb{R}^2$ as the \emph{bifurcation locus} of $\mathrm{Sp}(4)//\mathrm{Sp}(4)$. We define the bifurcation locus of $\mathrm{Sp}^\mathcal{I}(4)//\mathrm{GL}_2(\mathbb{R})$ in the same way; note that the map $\mathrm{Sp}^\mathcal{I}(4)//\mathrm{GL}_2(\mathbb{R})\rightarrow \mathrm{Sp}(4)//\mathrm{Sp}(4)$ maps the latter to the former.  

\begin{figure}
    \centering
    \includegraphics{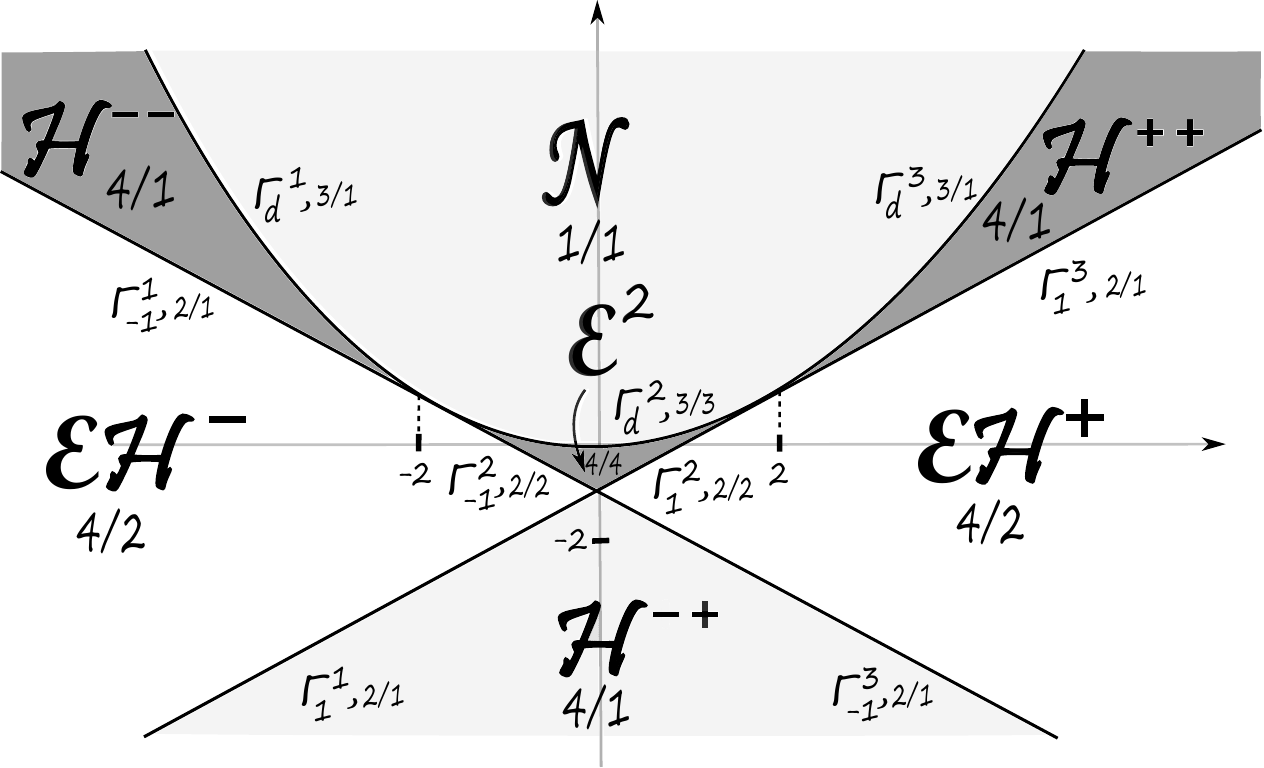}
    \caption{The branching locus and the seven components of $M_{2\times 2}(\mathbb{R})//GL_2(\mathbb{R})$. We denote the number of sheets on the interior of each component as $a/b$, where $a$ is the number of sheets of $\mathrm{Sp}^\mathcal{I}(4)//\mathrm{GL}_2(\mathbb{R})$, and $b$, that of $\mathrm{Sp}(4)//\mathrm{Sp}(4).$ Similarly, we indicate the number of sheets over the branch locus; cf.\ Figure \ref{fig:bifurcations}. This information gives a refinement of the ``Broucke stability diagram'' \cite{Broucke}.}
    \label{fig:branchlocus}
\end{figure}

The branch locus $\Gamma_d$ touches the
branching loci $\Gamma_1$ and $\Gamma_{-1}$ in the points
$$\Gamma_d \cap \Gamma_1=\big\{(2,1)\big\}, \quad \Gamma_d \cap \Gamma_{-1}=\big\{(-2,1)\big\}$$
at which the characteristic polynomial has a double root at $1$ respectively $-1$. Finally the
two branch loci $\Gamma_1$ and $\Gamma_{-1}$ intersect in the point
$$\Gamma_1 \cap \Gamma_{-1}=\big\{(0,-1)\big\}$$
at which the characteristic polynomial has a root at $1$ and $-1$, i.e., the matrix $A$ has
eigenvalues $1$ and $-1$. We abbreviate by
$$\Gamma=\Gamma_d\cup \Gamma_1 \cup \Gamma_{-1}$$
the full branch locus. Its complement decomposes into seven connected components
\begin{equation}\label{decom}
\mathbb{R}^2 \setminus \Gamma=\mathcal{E}^2 \cup \mathcal{EH}^+ \cup \mathcal{EH}^- \cup 
\mathcal{H}^{++} \cup \mathcal{H}^{-+} \cup \mathcal{H}^{--} \cup \mathcal{N}
\end{equation}
which we describe next; see Figure \ref{fig:branchlocus}.

The only bounded component in the decomposition the \emph{doubly elliptic component}
$$\mathcal{E}^2=\Big\{(\tau,\delta) \in \mathbb{R}^2: -2<\tau<2,\,\,\max\{-\tau-1,\tau-1\}<\delta<\tfrac{1}{4}
\tau^2\Big\}.$$
A matrix $A$ corresponding to this component has two distinct real eigenvalues
$-1<\mu_1<\mu_2<1$, while a matrix $M_{A,B,C}\in \mathrm{Sp}^\mathcal{I}(4)$ corresponding to
such a matrix $A$ has two pairs of eigenvalues on the unit circle
$(e^{i\theta_1},e^{-i\theta_1})$ and $(e^{i\theta_2},e^{-i\theta_2})$ with
$0<\theta_2<\theta_1<\pi$ which are related to the eigenvalues of $A$ in view of Lemma~\ref{cha} by
$$e^{i\theta_1}=\mu_1+i\sqrt{1-\mu_1^2}, \quad e^{i\theta_2}=\mu_2+i\sqrt{1-\mu_2^2}.$$
The \emph{elliptic/positive hyperbolic component} $\mathcal{EH}^+$ is given by
$$\mathcal{EH}^+=\Big\{(\tau,\delta) \in \mathbb{R}^2: \tau>0,\,\,\big|\delta+1\big|<\tau\Big\}.$$
A matrix $A$ for this component has two distinct real eigenvalues
$-1<\mu_1<1<\mu_2$, while a matrix $M_{A,B,C} \in \mathrm{Sp}^\mathcal{I}(4)$ corresponding
to $A$ has one pair of eigenvalues on the unit circle $(e^{i\theta}, e^{-i\theta})$ for
$0<\theta<\pi$ and a pair of positive real eigenvalues $\big(\lambda,\tfrac{1}{\lambda}\big)$
with $\lambda>1$ such that
$$e^{i\theta}=\mu_1+i\sqrt{1-\mu_1^2}, \quad \lambda=\mu_2+\sqrt{\mu_2^2-1}.$$
The \emph{elliptic/negative hyperbolic component} $\mathcal{EH}^-$ is given by
$$\mathcal{EH}^-=\Big\{(\tau,\delta) \in \mathbb{R}^2: \tau<0,\,\,\big|\delta+1\big|<-\tau\Big\}.$$
A matrix $A$ for this component has two distinct real eigenvalues
$\mu_1<-1<\mu_2<1$, while a matrix $M_{A,B,C} \in \mathrm{Sp}^\mathcal{I}(4)$ corresponding
to $A$ has one pair of eigenvalues on the unit circle $(e^{i\theta}, e^{-i\theta})$ for
$0<\theta<\pi$ and a pair of negative real eigenvalues $\big(\lambda,\tfrac{1}{\lambda}\big)$
with $\lambda<-1$ such that
$$e^{i\theta}=\mu_2+i\sqrt{1-\mu_2^2}, \quad \lambda=\mu_1-\sqrt{\mu_1^2-1}.$$
The \emph{negative/positive hyperbolic component} $\mathcal{H}^{-+}$ is given by
$$\mathcal{H}^{-+}=\Big\{(\tau,\delta) \in \mathbb{R}^2: \delta<-1,\,\,\big|\tau\big|<-\delta-1\Big\}.$$
A matrix $A$ for this component has two distinct real eigenvalues
$\mu_1<-1<1<\mu_2$, while a matrix $M_{A,B,C} \in \mathrm{Sp}^\mathcal{I}(4)$ corresponding
to $A$ has one pair negative real eigenvalues $\big(\lambda_1,\tfrac{1}{\lambda_1}\big)$
with $\lambda_1<-1$ and one pair of positive real eigenvalues
$\big(\lambda_2,\tfrac{1}{\lambda_2}\big)$ such that
$$\lambda_1=\mu_1-\sqrt{\mu_1^2-1}, \quad \lambda_2=\mu_2+\sqrt{\mu_2^2-1}.$$
The \emph{doubly positive hyperbolic component} $\mathcal{H}^{++}$ is given by
$$\mathcal{H}^{++}=\Big\{(\tau,\delta) \in \mathbb{R}^2: \tau>2,\,\,\tau-1<\delta<\tfrac{1}{4}
\tau^2\Big\}.$$
A matrix $A$ for this component has two distinct real eigenvalues
$1<\mu_1<\mu_2$, while a matrix $M_{A,B,C} \in \mathrm{Sp}^\mathcal{I}(4)$ corresponding
to $A$ has two pairs of positive real eigenvalues $\big(\lambda_1,\tfrac{1}{\lambda_1}\big)$
and $\big(\lambda_2,\tfrac{1}{\lambda_2}\big)$
with $1<\lambda_1<\lambda_2$ such that
$$\lambda_1=\mu_1+\sqrt{\mu_1^2-1}, \quad \lambda_2=\mu_2+\sqrt{\mu_2^2-1}.$$
The \emph{doubly negative hyperbolic component} $\mathcal{H}^{--}$ is given by
$$\mathcal{H}^{--}=\Big\{(\tau,\delta) \in \mathbb{R}^2: \tau<-2,\,\,-\tau-1<\delta<\tfrac{1}{4}
\tau^2\Big\}.$$
A matrix $A$ for this component has two distinct real eigenvalues
$\mu_1<\mu_2<-1$, while a matrix $M_{A,B,C} \in \mathrm{Sp}^\mathcal{I}(4)$ corresponding
to $A$ has two pairs of negative real eigenvalues $\big(\lambda_1,\tfrac{1}{\lambda_1}\big)$
and $\big(\lambda_2,\tfrac{1}{\lambda_2}\big)$
with $\lambda_1<\lambda_2<-1$ such that
$$\lambda_1=\mu_1-\sqrt{\mu_1^2-1}, \quad \lambda_2=\mu_2-\sqrt{\mu_2^2-1}.$$
Finally, the \emph{nonreal component} $\mathcal{N}$ is given by the region above the graph of the parabola
$\delta=\tfrac{1}{4}\tau^2$
$$\mathcal{N}=\Big\{(\tau,\delta) \in \mathbb{R}^2: \delta>\tfrac{1}{4}\tau^2\Big\}.$$
A matrix $A$ for this component has no real eigenvalues but a pair of two nonreal complex conjugated
eigenvalues $(\mu,\overline{\mu})$. A matrix $M_{A,B,C} \in \mathrm{Sp}^\mathcal{I}(4)$ has than a
quadruple of complex eigenvalues $\big(\lambda,\overline{\lambda},\tfrac{1}{\lambda},\tfrac{1}{\overline{\lambda}}\big)$ which are neither real nor lie on the unit circle where
$$\lambda=\mu+\sqrt{\mu^2-1}$$
where in this case $\sqrt{\mu^2-1}$ is the choice of a complex root of the complex number
$\mu^2-1$.
\\ \\
The union of the first six connected components in the decomposition (\ref{decom}) we abbreviate by
$$\mathcal{R}=\mathcal{E}^2 \cup \mathcal{EH}^+ \cup \mathcal{EH}^- \cup \mathcal{H}^{++} \cup
\mathcal{H}^{-+} \cup \mathcal{H}^{--}$$
and refer to it as the \emph{real part} of $\mathbb{R}^2 \setminus \Gamma$. With this notion we
have a decomposition
$$\mathbb{R}^2 \setminus \Gamma=\mathcal{R} \cup \mathcal{N}$$
into real and nonreal part. An equivalence class of matrices $[A]$ in the real part has two
distinct real eigenvalues, while in the nonreal part it has two nonreal eigenvalues which are
related to each other by complex conjugation. 
\\ \\
If $V \subset \mathbb{R}^2$ is an open subset we denote by
$\mathrm{Sp}(4)//\mathrm{Sp}(4)\big|_V$
the subset of $\mathrm{Sp}(4)//\mathrm{Sp}(4)$ consisting of all
$[[M_{A,B,C}]] \in \mathrm{Sp}(4)//\mathrm{Sp}(4)$ such that $[A]$ lies in $V$
and similarly $\mathrm{Sp}^\mathcal{I}(4)//\mathrm{GL}_2(\mathbb{R})\big|_V$. Outside
the branch locus $\Gamma$ the maps 
$$\mathrm{Sp}^\mathcal{I}(4)//\mathrm{GL}_2(\mathbb{R})\big|_{\mathbb{R}^2 \setminus \Gamma}
\to \mathrm{Sp}(4)//\mathrm{Sp}(4)|_{\mathbb{R}^2 \setminus \Gamma} \to
\mathbb{R}^2 \setminus \Gamma$$
are smooth coverings where the number of sheets however depends on the connected component in
$\mathbb{R}^2 \setminus \Gamma$. On this set it actually does not matter if one is considering
the GIT quotient or just the usual quotient.
\\ \\
Suppose now that $M_{A,B,C} \in \mathrm{Sp}^\mathcal{I}(4)$ and
$A$ has two distinct real eigenvalues, i.e., $[A] \in \mathcal{R}$.
Let $\mu$ be one of the eigenvalues of $A$. Its eigenspace $E_\mu \subset \mathbb{R}^2$ is 
then one-dimensional. In an eigenbasis the matrix
$A$ is diagonal and hence in view
of the equation $AB=BA$ it follows that $B$ leaves $E_\mu$ invariant. In particular,
there exists a real number $b_\mu$ such that for any eigenvector $v$ to the eigenvalue $\mu$
of $A$ we have
$$Bv=b_\mu v.$$
\begin{fed}\label{def:Bpos}
The eigenvalue $\mu$ of $A$ is called \emph{$B$-positive} if $b_\mu$ is positive and
\emph{$B$-negative} if $b_\mu$ is negative.
\end{fed}
Note that positivity and negativity of $b_\mu$ does not depend on
the choice of the eigenbasis, since $B$ transforms as a symmetric form
and therefore under change of the eigenbasis $b_\mu$ gets multiplied
by a positive number.
We now consider the elliptic case, i.e., $-1<\mu<1$. In particular, we must have
$[A] \in \mathcal{E}^2 \cup \mathcal{EH}^+ \cup \mathcal{EH}^-$. Then
$$\lambda=\mu+i\sqrt{1-\mu^2}$$
is an eigenvalue of the symplectic matrix $M_{A,B,C}$. 
\begin{lemma}\label{lemma:BposKrein}
In the elliptic case the eigenvalue $\mu$ of $A$ is $B$-positive (negative) if and only
if the eigenvalue $\lambda$ of $M_{A,B,C}$ is Krein-positive (negative).
\end{lemma}

We shall prove this lemma in Appendix \ref{app:Krein}. Note that it is crucial to take the positive sign for the imaginary part
of the eigenvalue $\lambda$. Its complex conjugate
$$\overline{\lambda}=\mu-i\sqrt{1-\mu^2}$$
is then another eigenvalue of the symplectic matrix $M_{A,B,C}$ of opposite Krein-type than
$\lambda$. That means that if $\mu$ is $B$-positive, then $\overline{\lambda}$ is
Krein-negative and if $\mu$ is $B$-negative, then $\overline{\lambda}$ is Krein-positive.
We further point out that the Krein-type of an eigenvalue of a symplectic matrix only depends on the conjugation class of the symplectic matrix. 
\\ \\
In the hyperbolic case where a real eigenvalue $\mu$ of $A$ satisfies $|\mu|>1$ there is 
no analogon of the Krein-type. On the other hand the $B$-type of an eigenvalue of a symplectic
matrix $M_{A,B,C} \in \mathrm{Sp}^\mathcal{I}(4)$ is defined in the hyperbolic case as well and
independent of the action of $\mathrm{GL}_2(\mathbb{R)}$ on $\mathrm{Sp}^\mathcal{I}(4)$. This
is the reason that in the hyperbolic case there are more sheets in the covering
$\mathrm{Sp}^\mathcal{I}(4)//\mathrm{GL}_2(\mathbb{R}) \to \mathbb{R}^2$, than in the branched covering
$\mathrm{Sp}(4)//\mathrm{Sp}(4)\to \mathbb{R}^2$.

\section{The characteristic polynomial}\label{sec:charpol} 

In this section we prove Lemma~\ref{cha}. For a matrix $M_{A,B,C}
\in \mathrm{Sp}^\mathcal{I}(2n)$ its characteristic polynomial is given by
$$
p_{A,B,C}(t)=\det\left(\begin{array}{cc}
A-tI & B\\
C & A^T-tI
\end{array}\right).
$$
Note that in view of (\ref{eq}) we have
\begin{eqnarray*}
\left(\begin{array}{cc}
A-tI & B\\
C & A^T-tI
\end{array}\right)\left(\begin{array}{cc}
A-tI & 0\\
-C & I
\end{array}\right)&=&\left(\begin{array}{cc}
(A-tI)^2-BC & B\\
CA-A^TC & A^T-tI
\end{array}\right)\\
&=&\left(\begin{array}{cc}
t^2 I-2t A+I & B\\
0 & A^T-tI
\end{array}\right).
\end{eqnarray*}
Taking determinants on both sides we obtain
\begin{eqnarray*}
p_{A,B,C}(t) \cdot \det(A-tI)&=&\det(t^2 I-2tA+I)\cdot \det(A^T-tI)\\
&=&\det(t^2 I-2tA+I)\cdot \det(A-tI)
\end{eqnarray*}
and hence
$$p_{A,B,C}(t)=\det(t^2I-2tA+I)=t^np_{-2A}\big(-t-\tfrac{1}{t}\big)$$
where $p_{-2A}$ is the characteristic polynomial of the matrix $-2A$. 
This proves Lemma~\ref{cha}.
\\ \\
Note that in the case $2n=4$ we have
$$p_{-2A}(s)=s^2-2(\mathrm{tr}A)s+4\det A$$
and therefore
$$p_{A,B,C}(t)=t^4-2(\mathrm{tr} A) t^3+2(1+2\det A)t^2-2(\mathrm{tr} A) t+1.$$

\section{Planar vs.\ spatial GIT quotients}

In this section, we explain the relationship between the GIT quotients for the two dimensional case (or planar case, i.e.\ $n=1$), and the four dimensional case (or spatial case, i.e\ $n=2$). Intuitively speaking, the spatial case behaves as a product of two planar cases (i.e.\ when two pairs of eigenvalues are independent of each other), except for the case where a non-real quadruple arises. Topologically, this means that the product of two copies of the GIT spaces for $n=1$ corresponds to the GIT space for $n=2$ with the non-real locus removed (although one needs to take a quotient by a $\mathbb{Z}_2$ action which forgets the order of the matrices which lie over the locus $\Gamma_d$). The details are as follows. 

Note that the product of two copies of the base of GIT sequence for $n=1$ is a copy of $\mathbb{R}^2$, and an element in this space corresponds to an ordered list of eigenvalues. For $n=2$, the base is parametrized by the trace and determinant of a $2\times2$-matrix. Therefore the map to consider is

$$F \colon \R^2 \to \R^2, (a,b) \mapsto (a+b,ab),$$

i.e.\ the map which associates to an ordered list of the two eigenvalues the trace and determinant of the matrix. In view of the inequality 
$$
(a+b)^2 \geq (a+b)^2-(a-b)^2=4ab,
$$
the image of $F$ precisely misses the non-real component $\mathcal{N}$ of the base of the GIT sequence for $n=2$. Moreover, on $\mathbb R^2$ we have the involution
$$
I \colon \mathbb R^2 \to \mathbb R^2, (a,b) \mapsto (b,a)
$$
interchanging the two eigenvalues. The map $F$ is invariant under the
involution $I$, i.e.\
$$
F \circ I=F.
$$
This reflects the fact that for the trace and determinant the order of
the eigenvalues does not play a role. The fixed point set of $I$ is the diagonal $\Delta$ in $\mathbb R^2$ which is mapped under $F$ precisely to the branching locus $\Gamma_d$. If $p_0:\mathrm{Sp}^{\mathcal{I}}(2n)//\mathrm{GL}_n(\mathbb R)\rightarrow\mathrm{Sp}(2n)//\mathrm{Sp}(2n)$ and  $p_1:\mathrm{Sp}(2n)//\mathrm{Sp}(2n)\rightarrow \mathrm{M}_{n\times n}(\mathbb R)//\mathrm{GL}_n(\mathbb{R})$ are the maps in the GIT sequence, and $p_2=p_1\circ p_0$, we conclude that
$$
(\mathrm{Sp}^{\mathcal{I}}(4)//\mathrm{GL}_2(\mathbb R))\backslash p_2^{-1}(\mathcal N)\simeq (\mathrm{Sp}^{\mathcal{I}}(2)//\mathrm{GL}_1(\mathbb R)) \times (\mathrm{Sp}^{\mathcal{I}}(2)//\mathrm{GL}_1(\mathbb R))/\sim, 
$$
$$
(\mathrm{Sp}(4)//\mathrm{Sp}(4))\backslash p_1^{-1}(\mathcal N)\simeq (\mathrm{Sp}(2)//\mathrm{Sp}(2))\times( \mathrm{Sp}(2)//\mathrm{Sp}(2))/\sim, 
$$
where the quotient identifies a pair $(M_1=M_{A_1,B_1,C_1},M_2=M_{A_2,B_2,C_2})$ with first blocks satisfying $A_1=A_2 \in \mathbb R$, with the pair $I(M_1,M_2)=(M_2,M_1)$.

\section{Higher-order bifurcations and pencils of lines}

In this section, for $n=2$, we consider the complete bifurcation locus, i.e.\ the locus of matrices having one eigenvalue which is a $k$-th root of unity for some $k$, and therefore corresponding to a $k$-fold bifurcation. We will see that this locus projects to the base of the GIT sequence as a line whose slope depends on $k$. More generally, we will consider the locus of matrices having a fixed eigenvalue. It turns out that the collection of such loci gives a pencil of lines in the plane, tangent to a parabola. The details are as follows.

\begin{figure}
    \centering
    \includegraphics{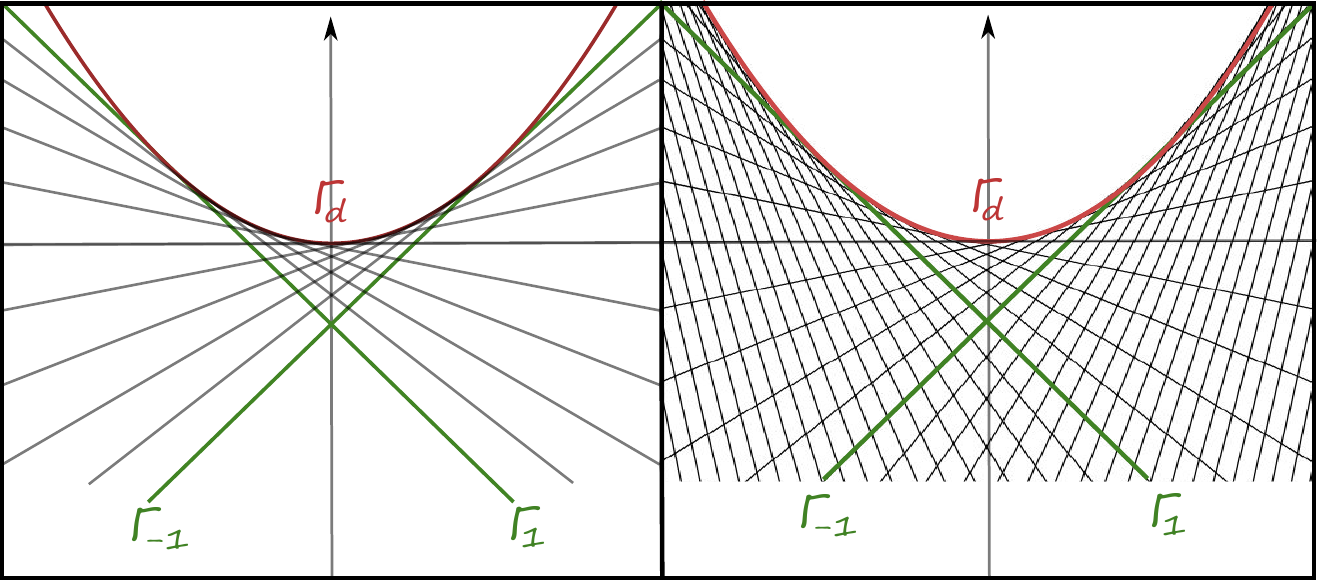}
    \caption{On the left, the \emph{elliptic} pencil of lines $\{\Gamma_\theta, \theta \in [0,2\pi)\}$. On the right, the complete pencil, also containing the \emph{hyperbolic} pencil $\{\Gamma_\lambda: \lambda\in \mathbb R\backslash[-1,1]\}$.}
    \label{fig:complete_bifurcation}
\end{figure}

Let $\lambda$ be an eigenvalue of $M=M_{A,B,C}$, which is a $k$-fold root of unity, i.e.\ it satisfies $\lambda^k=1$. By Lemma \ref{cha}, we have $a=a(\lambda)=\frac{1}{2}(\lambda+\frac{1}{\lambda})$ is an eigenvalue of $A$. If we write $\lambda=e^{2\pi i \ell/k}$, we have $a=\cos(2\pi \ell/k)$. Moreover, if $b$ is the remaining eigenvalue of $A$, its trace is $\tau=a+b$, and its determinant is $\delta=ab$, and we have the following equation
$$
\delta=a\tau-a^2=\cos(2\pi \ell/k)\tau-\cos(2\pi \ell/k)^2,
$$
which is a linear relation between $\delta$ and $\tau$, depending on $k$ and $\ell$. Note that the resulting lines in $\mathbb{R}^2$, denoted $\Gamma_{k,\ell}$, are all tangent to the parabola $\Gamma_d$, as no bifurcation can occur over the nonreal component $\mathcal{N}$. Moreover, two of such lines intersect at a point, consisting of the loci of curves whose two eigenvalues bifurcate with the orders corresponding to the lines. More generally, one can consider the locus $\Gamma_\theta$ consisting of matrices having $e^{2\pi i\theta}$ as an eigenvalue. The same computation as above shows that this is the line
$$
\delta=\cos(2\pi \theta)\tau-\cos(2\pi \theta)^2,
$$
tangent to $\Gamma_d$ at $\tau=2\cos(2\pi \theta)$.
The collection $\{\Gamma_\theta:\theta \in [0,2\pi)\}$ is a pencil of lines, tangent to the parabola $\Gamma_d$, but with slopes varying only in $[-1,1]$. The intersection of $\Gamma_\theta$ with $\Gamma_{\theta'}$ lies in $\mathcal{E}^2$, and is the locus of matrices having eigenvalues $e^{\pm 2\pi i\theta},e^{\pm 2\pi i\theta'}$. Analogously, we can consider the case where $\lambda$ is real, i.e.\ hyperbolic, in which case the resulting line $\Gamma_\lambda=\{\delta=a(\lambda)\tau-a(\lambda)^2\}$, consisting of the locus of matrices with eigenvalue $\lambda$, is tangent to $\Gamma_d$ at $\tau=2a$; note that $\Gamma_\lambda=\Gamma_{1/\lambda}$. The slope $a(\lambda)$ is greater than 1 (resp.\ smaller than $-1$) if and only if $\lambda$ is positive (resp.\ negative) hyperbolic. The intersections between any of the lines $\Gamma_\theta,\Gamma_\lambda$ again have the obvious interpretation. See Figure \ref{fig:complete_bifurcation}.

\section{Normal forms}

In this section we describe normal forms for matrices $M_{A,B,C} \in \mathrm{Sp}^\mathcal{I}(4)$
under the action of $\mathrm{GL}_2(\mathbb{R})$ given by (\ref{act}). Our normal forms still
lie in $\mathrm{Sp}^\mathcal{I}(4)$, therefore with the exception of the doubly elliptic case
they differ rather from standard normal forms 
of symplectic matrices as for example explained in \cite{long}.

\subsection{The regular cases}

In this section we assume that our symplectic matrix does not lie
over the branch locus. Under this assumption the orbits are closed
and there is no difference between the equivalence class in the GIT quotient
and the usual quotient. 
\\ \\
We first consider the case where $A$ has two different real eigenvalues, i.e., $[A] \in \mathcal{R}$.
In this case $A$ is diagonalizable and after acting with $\mathrm{GL}_2(\mathbb{R})$ we can
assume without loss of generality that $A$ is actually diagonal
$$A=\left(\begin{array}{cc}
\mu_1 & 0\\
0 & \mu_2
\end{array}\right)$$
while we order the eigenvalues in increasing order $\mu_1<\mu_2$. In particular, we have
that $A=A^T$ and therefore
$$AB=BA, \quad AC=CA,$$
i.e., $A$ commutes with both matrices $B$ and $C$. Since the two eigenvalues of $A$ are different,
this implies that $B$ and $C$ are diagonal as well
$$B=\left(\begin{array}{cc}
b_1 & 0\\
0 & b_2
\end{array}\right), \quad C=\left(\begin{array}{cc}
c_1 & 0\\
0 & c_2
\end{array}\right).$$
From the equation $A^2-BC=I$ we obtain
$$b_1c_1=\mu_1^2-1, \quad b_2c_2=\mu_2^2-1.$$
After acting by the matrix
$$R=\left(\begin{array}{cc}
\big|\frac{b_1}{c_1}\big|^{1/4} & 0\\
0 & \big|\frac{b_2}{c_2}\big|^{1/4}
\end{array}\right)$$
according to (\ref{act}) we can achieve in addition that
$$|b_1|=|c_1|, \quad |b_2|=|c_2|$$
so that we have
$$b_1^2=|\mu_1^2-1|, \quad b_2^2=|\mu_2^2-1|.$$
In the following we discuss the six connected components of the real part $\mathcal{R}$ one by one.
We start by assuming that $[A]$ lies in the double elliptic component $\mathcal{E}^2$. In this case
the two real eigenvalues of $A$ satisfy $-1<\mu_1<\mu_2<1$. In particular, there exist unique
angles $\theta_1,\theta_2 \in (0,\pi)$ such that
$$\mu_1=\cos \theta_1, \quad \mu_2=\cos \theta_2.$$
We then have
$$|b_1|=|c_1|=\sqrt{1-\cos^2 \theta_1}=|\sin \theta_1|, \quad
 |b_2|=|c_2|=\sqrt{1-\cos^2 \theta_2}=|\sin \theta_2|$$
and the signs of $b_1$ and $c_1$ as well as those of $b_2$ and $c_2$ have to be opposite. 
We therefore get the following four normal forms
$$\left(\begin{array}{cccc}
\cos \theta_1 &      0       & -\sin \theta_1 & 0\\
0             & \cos\theta_2 &            0   & -\sin \theta_2\\
\sin \theta_1 &     0        & \cos \theta_1  & 0\\
0             & \sin\theta_2 &            0   & \cos \theta_2
\end{array}\right), 
\left(\begin{array}{cccc}
\cos \theta_1 &      0       & \sin \theta_1 & 0\\
0             & \cos\theta_2 &            0   & -\sin \theta_2\\
-\sin \theta_1 &     0        & \cos \theta_1  & 0\\
0             & \sin\theta_2 &            0   & \cos \theta_2
\end{array}\right),$$ 
$$\left(\begin{array}{cccc}
\cos \theta_1 &      0       & -\sin \theta_1 & 0\\
0             & \cos\theta_2 &            0   & \sin \theta_2\\
\sin \theta_1 &     0        & \cos \theta_1  & 0\\
0             & -\sin\theta_2 &            0   & \cos \theta_2
\end{array}\right), 
\left(\begin{array}{cccc}
\cos \theta_1 &      0       & \sin \theta_1 & 0\\
0             & \cos\theta_2 &            0   & \sin \theta_2\\
-\sin \theta_1 &     0        & \cos \theta_1  & 0\\
0             & -\sin\theta_2 &            0   & \cos \theta_2
\end{array}\right).$$
We therefore see that the fiber of $\mathrm{Sp}^\mathcal{I}(4)//\mathrm{GL}_2(\mathbb{R})$ over $[A]$
consists of four points. Moreover, in view of Lemma~\ref{lemma:BposKrein} all these matrices are distinguished symplectically
by the Krein-type of their eigenvalues. Therefore the fiber of $\mathrm{Sp}(4)//\mathrm{Sp}(4)$
consists of four points as well. 
\\ \\
We next discuss the case that $[A]$ lies in the 
elliptic/positive hyperbolic component $\mathcal{EH}^+$. In this
case the two real eigenvalues of $A$ satisfy $-1<\mu_1<1<\mu_2$.
Hence there exist unique $\theta_1 \in (0,\pi)$ and
$\theta_2 \in (0,\infty)$ such that
$$\mu_1=\cos \theta_1, \quad \mu_2=\cosh \theta_2.$$
Hence
$$|b_1|=|c_1|=\sqrt{1-\cos^2 \theta_1}=|\sin \theta_1|,
\quad |b_2|=|c_2|=\sqrt{\cosh^2 \theta_2-1}=|\sinh \theta_2|.$$
Where the signs of $b_1$ and $c_1$ are opposite the signs of $b_2$ and $c_2$ agree. We obtain the
following four normal forms.
$$\left(\begin{array}{cccc}
\cos \theta_1 &      0       & -\sin \theta_1 & 0\\
0             & \cosh\theta_2 &            0   & -\sinh \theta_2\\
\sin \theta_1 &     0        & \cos \theta_1  & 0\\
0             & -\sinh\theta_2 &            0   & \cosh \theta_2
\end{array}\right), 
\left(\begin{array}{cccc}
\cos \theta_1 &      0       & -\sin \theta_1 & 0\\
0             & \cosh\theta_2 &            0   & \sinh \theta_2\\
\sin \theta_1 &     0        & \cosh \theta_1  & 0\\
0             & \sinh\theta_2 &            0   & \cosh \theta_2
\end{array}\right),$$ 
$$\left(\begin{array}{cccc}
\cos \theta_1 &      0       & \sin \theta_1 & 0\\
0             & \cosh\theta_2 &            0   & -\sinh \theta_2\\
-\sin \theta_1 &     0        & \cos \theta_1  & 0\\
0             & -\sinh\theta_2 &            0   & \cosh \theta_2
\end{array}\right), 
\left(\begin{array}{cccc}
\cos \theta_1 &      0       & \sin \theta_1 & 0\\
0             & \cosh\theta_2 &            0   & \sinh \theta_2\\
-\sin \theta_1 &     0        & \cos \theta_1  & 0\\
0             & \sinh\theta_2 &            0   & \cosh \theta_2
\end{array}\right).$$
The fiber of $\mathrm{Sp}^\mathcal{I}(4)//\mathrm{GL}_2(\mathbb{R})$ over $[A]$
consists again of four points. However, the two matrices on
the first line are symplectically conjugated via the symplectic matrix
$$\left(\begin{array}{cccc}
1 & 0 & 0 & 0\\
0 & 0 & 0 & -1\\
0 & 0 & 1 & 0\\
0 & 1 & 0 & 0
\end{array}\right)$$
and the same holds for the two matrices on the second line. On the other hand
the matrices on the first line and the ones on the second line are symplectically still
distinguished by the Krein-type of the first eigenvalue $\mu_1=\cos \theta_1$. Hence the fiber of $\mathrm{Sp}(4)//\mathrm{Sp}(4)$ consists
over $[A]$ consists of two points. 
\\ \\
In the elliptic/negative hyperbolic case the eigenvalues satisfy 
$\mu_1<-1<\mu_2<1$ so that we write
$$\mu_1=-\cosh \theta_1, \quad \mu_2=\cos \theta_2$$
for $\theta_1 \in (0,\infty)$ and $\theta_2 \in (0,\pi)$. We have the following four normal
forms
$$\left(\begin{array}{cccc}
-\cosh \theta_1 &      0       & -\sinh \theta_1 & 0\\
0             & \cos\theta_2 &            0   & -\sin \theta_2\\
-\sinh \theta_1 &     0        & -\cosh \theta_1  & 0\\
0             & \sin\theta_2 &            0   & \cos \theta_2
\end{array}\right), 
\left(\begin{array}{cccc}
-\cosh \theta_1 &      0       & \sinh \theta_1 & 0\\
0             & \cos\theta_2 &            0   & -\sin \theta_2\\
\sinh \theta_1 &     0        & -\cosh \theta_1  & 0\\
0             & \sin\theta_2 &            0   & \cos \theta_2
\end{array}\right),$$ 
$$\left(\begin{array}{cccc}
-\cosh \theta_1 &      0       & -\sinh \theta_1 & 0\\
0             & \cos\theta_2 &            0   & \sin \theta_2\\
-\sinh \theta_1 &     0        & -\cosh \theta_1  & 0\\
0             & -\sin\theta_2 &            0   & \cos \theta_2
\end{array}\right), 
\left(\begin{array}{cccc}
-\cosh \theta_1 &      0       & \sinh \theta_1 & 0\\
0             & \cos\theta_2 &            0   & \sin \theta_2\\
\sinh \theta_1 &     0        & -\cosh \theta_1  & 0\\
0             & -\sin\theta_2 &            0   & \cos \theta_2
\end{array}\right).$$
Again on each line the matrices are symplectically conjugated while on different lines
they are symplectically distinguished by the Krein-type of the second eigenvalue. Therefore
the fiber of $\mathrm{Sp}^\mathcal{I}(4)//\mathrm{GL}_2(\mathbb{R})$ over $[A]$ has again
four points while the one of $\mathrm{Sp}(4)//\mathrm{Sp}(4)$ has just two points. 
\\ \\
In the doubly positive hyperbolic case $\mathcal{H}^{++}$ we obtain the following four normal forms
$$\left(\begin{array}{cccc}
\cosh \theta_1 &      0       & -\sinh \theta_1 & 0\\
0             & \cosh\theta_2 &            0   & -\sinh \theta_2\\
-\sinh \theta_1 &     0        & \cosh \theta_1  & 0\\
0             & -\sinh\theta_2 &            0   & \cosh \theta_2
\end{array}\right), 
\left(\begin{array}{cccc}
\cosh \theta_1 &      0       & \sinh \theta_1 & 0\\
0             & \cosh\theta_2 &            0   & -\sinh \theta_2\\
\sinh \theta_1 &     0        & \cosh \theta_1  & 0\\
0             & -\sinh\theta_2 &            0   & \cosh \theta_2
\end{array}\right),$$ 
$$\left(\begin{array}{cccc}
\cosh \theta_1 &      0       & -\sinh \theta_1 & 0\\
0             & \cosh\theta_2 &            0   & \sinh \theta_2\\
-\sinh \theta_1 &     0        & \cosh \theta_1  & 0\\
0             & \sinh\theta_2 &            0   & \cosh \theta_2
\end{array}\right), 
\left(\begin{array}{cccc}
\cosh \theta_1 &      0       & \sinh \theta_1 & 0\\
0             & \cosh\theta_2 &            0   & \sinh \theta_2\\
\sinh \theta_1 &     0        & \cosh \theta_1  & 0\\
0             & \sinh\theta_2 &            0   & \cosh \theta_2
\end{array}\right).$$
which now are all symplectically conjugated. 
Similarly are the negative/positive hyperbolic case $\mathcal{H}^{-+}$ and the doubly negative hyperbolic case $\mathcal{H}^{--}$. 
The only difference is that in the negative/positive hyperbolic
case $\cosh \theta_1$ has to be replaced by $-\cosh \theta_1$ and
in the doubly negative hyperbolic case both 
$\cosh \theta_1$ and $\cosh \theta_2$ get a minus sign. 
\\ \\
A different treatment needs the nonreal case, where $[A] \in \mathcal{N}$. In this case
$A$ has two nonreal eigenvalues which are complex conjugate of each other
$$\mu=re^{i\theta}, \quad \overline{\mu}=re^{-i\theta}$$
for $r>0$ and $\theta \in (0,\pi)$. After conjugation we can assume that the matrix
$A$ is the composition of a dilation by $r$ and a rotation by $\theta$
$$A=\left(\begin{array}{cc}
r \cos \theta & -r \sin \theta\\
r \sin \theta & r \cos \theta
\end{array} \right).$$
The matrix $B$ is symmetric and transforms as a bilinear form. Hence after a further rotation
which does not affect the matrix $A$ we can assume that $B$ is diagonal
$$B=\left(\begin{array}{cc}
b_1 & 0\\
0   & b_2
\end{array}\right).$$
The equation $AB=BA^T$ implies that $b_2=-b_1$ and after a further
dilation which still does not affect $A$ we can assume that
$B$ has the form
$$B=\left(\begin{array}{cc}
1 & 0\\
0 & -1
\end{array}\right),$$
i.e., $B$ is just an orthogonal reflection at the first axis. The equation $A^2-BC=I$ implies that
$$C=B(A^2-I)=\left(\begin{array}{cc}
r^2 \cos 2\theta -1 & -r^2\sin 2\theta\\
-r^2 \sin 2\theta & r^2 \cos 2\theta -1
\end{array}\right).$$
Therefore we obtain the unique canonical form
$$\left(\begin{array}{cccc}
r \cos \theta & -r \sin \theta & 1 & 0\\
r\sin \theta & r\cos \theta     & 0 & -1\\
r^2 \cos 2\theta-1 & -r^2\sin 2\theta & r \cos \theta & r\sin \theta\\
-r^2 \sin 2\theta  & r^2 \cos 2\theta-1 & -r\sin \theta & r \cos \theta
\end{array}\right).$$
In particular, since the canonical form is unique we see that both fibers
of $\mathrm{Sp}^\mathcal{I}(4)//\mathrm{GL}_2(\mathbb{R})$ and $\mathrm{Sp}(4)//\mathrm{Sp}(4)$
over $[A]$ consist of a single point, that means over the nonreal component the coverings are
just trivial, i.e., homeomorphisms. 

\subsection{The branch locus}

In this section we discuss normal forms over the branch locus. Over the branch locus not all orbits are closed and there is can be a difference between normal forms for the GIT quotient and the usual quotient. The branch
locus itself has three singular points at $(2,1)$, $(-2,1)$, and
$(0,-1)$. On the complement of the singularities the branch locus
consists of nine connected components all homeomorphic to an open interval. 
\\ \\
Recall that we abbreviated the first branch locus by $\Gamma_d$ which
is given by $\Gamma_d=\big\{\big(\tau,\tfrac{1}{4}\tau^2\big):\tau \in
\mathbb{R}\big\}$. It contains the two singularities $(2,1)$ and $(2,-1)$.
Its complement decomposes into three connected components
$$\Gamma_d\setminus \big\{(2,1) \cup (-2,1)\big\}
=\Gamma_d^1 \cup \Gamma_d^2 \cup \Gamma_d^3,$$
where 
\begin{eqnarray*}
\Gamma_d^1&=&\Big\{\big(\tau,\tfrac{1}{4}\tau^2\big):\tau<-2\Big\},\\
\Gamma_d^2&=&\Big\{\big(\tau,\tfrac{1}{4}\tau^2\big):-2<\tau<2\Big\},\\
\Gamma_d^3&=&\Big\{\big(\tau,\tfrac{1}{4}\tau^2\big):\tau>2\Big\}.
\end{eqnarray*}
Similarly we have the decomposition
$$\Gamma_1 \setminus \big\{(2,1),(0,-1)\big\}
=\Gamma_1^1 \cup \Gamma_1^2 \cup \Gamma_1^3$$
with
\begin{eqnarray*}
\Gamma_1^1&=&\Big\{(\tau,\tau-1): \tau<0\Big\}\\
\Gamma_1^2&=&\Big\{(\tau,\tau-1): 0<\tau<2\Big\}\\
\Gamma_1^3&=&\Big\{(\tau,\tau-1): \tau>2\Big\}
\end{eqnarray*}
as well as
$$\Gamma_{-1} \setminus \big\{(-2,1),(0,-1)\big\}
=\Gamma_{-1}^1 \cup \Gamma_{-1}^2 \cup \Gamma_{-1}^3$$
with
\begin{eqnarray*}
\Gamma_{-1}^1&=&\Big\{(\tau,-\tau-1): \tau<-2\Big\}\\
\Gamma_{-1}^2&=&\Big\{(\tau,-\tau-1): -2<\tau<0\Big\}\\
\Gamma_{-1}^3&=&\Big\{(\tau,-\tau-1): \tau>0\Big\}.
\end{eqnarray*}
With these notions the nonsingular part of the branch locus decomposes into connected components as follows
$$\Gamma \setminus \big\{(2,1),(-2,1),(0,-1)\big\}
=\Gamma_d^1 \cup \Gamma_d^2 \cup \Gamma_d^3\cup\Gamma_1^1 \cup \Gamma_1^2 \cup \Gamma_1^3\cup \Gamma_{-1}^1 \cup \Gamma_{-1}^2 \cup \Gamma_{-1}^3.$$
We first discuss the normal forms over the nonsingular part of the branch
locus. Hence we assume that $M_{A,B,C} \in \mathrm{Sp}^\mathcal{I}(4)$
with $[A] \in \Gamma \setminus \{(2,1),(-2,1),(0,-1)\}$.
We first assume that 
$[A] \in \Gamma_d \setminus \{(2,1),(-2,-1)\}$, i.e., $A$ has one
real eigenvalue $\mu$ different from $\pm 1$ of algebraic multiplicity 
two. The geometric multiplicity of the eigenvalue $\mu$ is one or two.
We first explain that by going over to the GIT quotient we can assume
that its geometric multiplicity is two as well. To see that suppose
that the geometric multiplicity of $\mu$ is one. 
After acting with $\mathrm{GL}_2(\mathbb{R})$ we can assume that $A$
is a Jordan block
$$A=\left(\begin{array}{cc}
\mu & 1\\
0 & \mu
\end{array}\right).$$
From the equation $AB=BA^T$ we infer that the symmetric matrix $B$ has the form
$$B=\left(\begin{array}{cc}
b_2 & b_1\\
b_1 & 0\end{array}\right)$$
where from the equation $A^TC=CA$ we deduce that $C$ has the form
$$C=\left(\begin{array}{cc}
0 & c_1\\
c_1 & c_2\end{array}\right).$$
For $\epsilon>0$ we consider the matrix
$$R_\epsilon=\left(\begin{array}{cc}
\epsilon & 0\\
0& \frac{1}{\epsilon}
\end{array}\right).$$
We have
\begin{eqnarray*}
R_\epsilon A R_\epsilon^{-1}&=& \left(\begin{array}{cc}
\mu & \epsilon^2\\
0 & \mu \end{array}\right)\\
R_\epsilon B R_\epsilon^T&=&\left(\begin{array}{cc}
\epsilon^2 b_2 & b_1\\
b_1 & 0\end{array}\right)\\
(R_\epsilon^T)^{-1} C R_\epsilon^{-1}&=&\left(\begin{array}{cc}
0 & c_1\\
c_1 & \epsilon^2 c_2\end{array}\right)
\end{eqnarray*}
and therefore
\begin{eqnarray*}
\lim_{\epsilon \to 0}R_\epsilon A R_\epsilon^{-1}&=& \left(\begin{array}{cc}
\mu & 0\\
0 & \mu \end{array}\right)\\
\lim_{\epsilon \to 0}R_\epsilon B R_\epsilon^T&=&\left(\begin{array}{cc}
0 & b_1\\
b_1 & 0\end{array}\right)\\
\lim_{\epsilon \to 0}(R_\epsilon^{T})^{-1} C R_\epsilon^{-1}&=&\left(\begin{array}{cc}
0 & c_1\\
c_1 & 0\end{array}\right).
\end{eqnarray*}
This shows that by going over to the GIT-quotient we can assume without
loss of generality that the geometric multiplicity of the eigenvalue 
$\mu$ is two as well. In this case the matrix $A$ is diagonal
$$A=\left(\begin{array}{cc}
\mu & 0\\
0 & \mu
\end{array}\right).$$
In fact $A$ is just a scalar multiple of the identity matrix and therefore
a fixed point of the action of $\mathrm{GL}_2(\mathbb{R})$ by conjugation.
The matrix $B$ transforms as a bilinear form and since any bilinear form
can be diagonalized we can assume after conjugation that $B$ is also diagonal. Since $\mu \neq \pm 1$ the formula $A^2-BC=I$ implies that $B$ has
to be nonsingular and therefore $C$ has to be diagonal as well. 
The discussion of normal forms is now analogous
to the real regular case. Different from the real regular case there are
only three normal forms and not four. These are in one to one correspondence with the signature of $B$. If the signature of $B$ is one in the real regular case there were still two different normal forms which were 
distinguished on which eigenspace of $A$ the matrix $B$ was positive and on
which it was negative. Since the two eigenvalues now coincide this distinction is not possible anymore. 
\\ \\
If $[A]$ lies in $\Gamma_d^2$, i.e., in the intersection of the
closure of the double elliptic component and the nonreal component, 
we have the following three normal forms for the eigenvalue
$\mu=\cos \theta$ with $\theta \in (0,\pi)$
$$\left(\begin{array}{cccc}
\cos \theta &      0       & -\sin \theta& 0\\
0             & \cos\theta &            0   & -\sin \theta\\
\sin \theta&     0        & \cos \theta  & 0\\
0             & \sin\theta&            0   & \cos \theta
\end{array}\right), 
\left(\begin{array}{cccc}
\cos \theta &      0       & \sin \theta & 0\\
0             & \cos\theta &            0   & -\sin \theta\\
-\sin \theta &     0        & \cos \theta & 0\\
0             & \sin\theta &            0   & \cos \theta
\end{array}\right),$$ 
$$
\left(\begin{array}{cccc}
\cos \theta&      0       & \sin \theta & 0\\
0             & \cos\theta &            0   & \sin \theta\\
-\sin \theta &     0        & \cos \theta  & 0\\
0             & -\sin\theta &            0   & \cos \theta
\end{array}\right).$$
Hence over $\Gamma^2_d$ the branch cover $\mathrm{Sp}^\mathcal{I}(4)//
\mathrm{GL}_2(\mathbb{R})$ has three branches. Moreover, the three normal
forms are distinguished symplectically by their Krein-type so that
$\mathrm{Sp}(4)//\mathrm{Sp}(4)$ over $\Gamma^2_d$ has three branches as well.
\\ \\
If $[A]$ lies in $\Gamma^3_d$, i.e., the intersection of the closures of the doubly positive hyperbolic component and the nonreal component we have the following three normal forms for $\mu=\cosh \theta$ with $\theta \in (0,\infty)$
$$\left(\begin{array}{cccc}
\cosh \theta &      0       & -\sinh \theta & 0\\
0             & \cosh\theta &            0   & -\sinh \theta\\
-\sinh \theta &     0        & \cosh \theta & 0\\
0             & -\sinh\theta &            0   & \cosh \theta
\end{array}\right), 
\left(\begin{array}{cccc}
\cosh \theta&      0       & \sinh \theta & 0\\
0             & \cosh\theta &            0   & -\sinh \theta\\
\sinh \theta &     0        & \cosh \theta  & 0\\
0             & -\sinh\theta &            0   & \cosh \theta
\end{array}\right),$$ 
$$
\left(\begin{array}{cccc}
\cosh \theta &      0       & \sinh \theta & 0\\
0             & \cosh\theta &            0   & \sinh \theta\\
\sinh \theta &     0        & \cosh \theta_1  & 0\\
0             & \sinh\theta &            0   & \cosh \theta
\end{array}\right).$$
In particular, over $\Gamma^3_d$ the branched cover
$\mathrm{Sp}^\mathcal{I}(4)//
\mathrm{GL}_2(\mathbb{R})$ has three branches. On the other hand the above
three normal forms are symplectically conjugated and therefore
$\mathrm{Sp}(4)//\mathrm{Sp}(4)$ has just one branch. A similar picture 
happens over $\Gamma^1_d$, i.e., the intersection of the closures
of the negative hyperbolic component and the nonreal component. There one
just needs to replace $\cosh \theta$ by $-\cosh \theta$ in the previous discussion.
\\ \\
We next assume that $[A] \in \Gamma_1 \setminus \{(2,1),(0,-1)\}$, i.e.,
$A$ has one eigenvalue equal to $1$ and another real eigenvalue $\mu \neq \pm 1$. In particular, $A$ is diagonalizable and after conjugation we can assume that $A$ has the form
$$A=\left(\begin{array}{cc}
  1   &  0\\
0     & \mu
\end{array}\right).$$
In particular, we have $A=A^T$ and therefore $B$ and $C$ commute with
$A$. This implies that they are diagonal as well
$$B=\left(\begin{array}{cc}
 b_1    &  0\\
0     & b_2 
\end{array}\right), \quad 
C=\left(\begin{array}{cc}
 c_1    &  0\\
0     & c_2 
\end{array}\right).$$
In view of $A^2-BC=I$ we obtain
$$b_1 c_1=0, \qquad b_2c_2=\mu^2-1.$$
The first equation implies that $b_1$ or $c_1$ is zero. We next explain that
by going over to the GIT quotient we can assume that both $b_1$ and $c_1$
are zero. To see that we first assume that $c_1=0$ but $b_1 \neq 0$ so that we have
$$A=\left(\begin{array}{cc}
  1   &  0\\
  0   & \mu
\end{array}\right), \quad B=\left(\begin{array}{cc}
 b_1    &  0\\
0     & b_2 
\end{array}\right), \quad 
C=\left(\begin{array}{cc}
 0    &  0\\
0     & c_2 
\end{array}\right).$$
For $\epsilon>0$ we consider the family of matrices
$$R_\epsilon=\left(\begin{array}{cc}
  \epsilon   &  0\\
   0  & 1
\end{array}\right).$$
Acting with $R_\epsilon$ on the triple of matrices above we obtain
\begin{eqnarray*}
R_\epsilon A R_\epsilon^{-1}&=& A\\
R_\epsilon B R_\epsilon^T&=&\left(\begin{array}{cc}
\epsilon^2 b_1 & 0\\
0 & b_2\end{array}\right)\\
(R_\epsilon^T)^{-1} C R_\epsilon^{-1}&=&C
\end{eqnarray*}
and therefore
\begin{eqnarray*}
\lim_{\epsilon \to 0}R_\epsilon A R_\epsilon^{-1}&=& A\\
\lim_{\epsilon \to 0}R_\epsilon B R_\epsilon^T&=&\left(\begin{array}{cc}
0 & 0\\
0 & b_2\end{array}\right)\\
\lim_{\epsilon \to 0}(R_\epsilon^{T})^{-1} C R_\epsilon^{-1}&=&C.
\end{eqnarray*}
This shows that we can assume that $b_1=0$. Similarly we see that
we can assume as well that $c_1=0$, by using instead in the above argument
the family of matrices
$$R_\epsilon=\left(\begin{array}{cc}
  \frac{1}{\epsilon}   &  0\\
   0  & 1
\end{array}\right).$$
After using the action of $\mathrm{GL}_2(\mathbb{R})$ once more we can additionally
assume that
$$|b_2|=|c_2|=\sqrt{|\mu^2-1|}.$$
If $[A] \in \Gamma_1^2$, i.e., if $[A]$ lies in the intersection of the
closures of the doubly elliptic component and the elliptic/positive hyperbolic component, we have the following two normal forms
for $\mu=\cos \theta$ with $\theta \in (0,\pi)$
$$\left(\begin{array}{cccc}
1 &      0       & 0 & 0\\
0             & \cos\theta &            0   & -\sin \theta\\
0 &     0        & 1  & 0\\
0             & \sin\theta&            0   & \cos \theta
\end{array}\right), 
\left(\begin{array}{cccc}
1 &      0       & 0 & 0\\
0             & \cos\theta &            0   & \sin \theta\\
0 &     0        & 1 & 0\\
0             & -\sin\theta &            0   & \cos \theta
\end{array}\right).$$ 
In particular, over $\Gamma_1^2$ the branched cover
$\mathrm{Sp}^\mathcal{I}(4)//\mathrm{GL}_2(\mathbb{R})$ has
two branches. This two branches are distinguished by their Krein-type and
therefore not symplectically conjugated, so that 
$\mathrm{Sp}(4)//\mathrm{Sp}(4)$ has two branches as well.
\\ \\
If $[A] \in \Gamma_1^3$, i.e., if $[A]$ lies in the intersection of the
closures of the doubly positive hyperbolic component and the elliptic/positive hyperbolic
component, the two normal forms are for $\mu=\cosh \theta$ with
$\theta \in (0,\infty)$
$$\left(\begin{array}{cccc}
1 &      0       & 0 & 0\\
0             & \cosh\theta &            0   & -\sinh \theta\\
0 &     0        & 1  & 0\\
0             & -\sinh\theta&            0   & \cosh \theta
\end{array}\right), 
\left(\begin{array}{cccc}
1 &      0       & 0 & 0\\
0             & \cosh\theta &            0   & \sinh \theta\\
0 &     0        & 1 & 0\\
0             & \sinh\theta &            0   & \cosh \theta
\end{array}\right).$$ 
Again $\mathrm{Sp}^\mathcal{I}(4)/\mathrm{GL}_2(\mathbb{R})$ has two
branches over $\Gamma_1^3$ but now the two branches are symplectically
conjugated and there is just one branch of
$\mathrm{Sp}(4)//\mathrm{Sp}(4)$. The case where 
$[A] \in \Gamma_1^1$, i.e., where $[A]$ lies in the intersection of
the closures of the elliptic/negative hyperbolic component and the
negative/positive hyperbolic component is similar, one just needs to replace
$\cosh \theta$ by $-\cosh \theta$.
\\ \\
Finally the discussion where $[A] \in \Gamma_{-1}\setminus \{(-2,1),(0,-1)\}$ is analogous to the one where
$[A] \in \Gamma_1\setminus \{(2,1),(0,-1)\}$. The only difference is that
one has to replace $1$ by $-1$. This finishes the description of the
branched covers over the nonsingular part of the branch locus. 
\\ \\
It remains to consider the singular part of the branch locus namely the
three points $(2,1)$, $(-2,1)$, and $(0,-1)$. We start with the point
$(2,1)$. In this case $A$ has only $1$ as eigenvalue with algebraic multiplicity two. If the geometric multiplicity is two as well, then $A$ is
the identity matrix. If the geometric multiplicity is one, then $A$ is
conjugated to the $2 \times 2$-Jordan block with $1$ on the diagonal. We explain that in either case the $4\times 4$-identity matrix $I_4$
lies in the closure of the $\mathrm{GL}_2(\mathbb{R})$-orbit
of $M_{A,B,C}$. For that purpose we first consider the case where $A$ is
the Jordan block
$$A=\left(\begin{array}{cc}
   1  & 1 \\
0     & 1
\end{array}\right).$$
From the equations $AB=BA^T$ and $A^TC=CA$ we infer that the symmetric
matrices $B$ and $C$ simplify to
$$B=\left(\begin{array}{cc}
    b_2 & b_1 \\
    b_1 & 0
\end{array}\right), \quad C=\left(\begin{array}{cc}
    0 & c_1 \\
    c_1 & c_2
\end{array}\right).$$
From the equation $A^2-BC=I$ we obtain
$$b_1c_1=0,\qquad c_1b_2+c_2b_1=2.$$
From the first equation we see that $b_1$ or $c_1$ vanishes. 
We first consider the case where $b_1=0$, so that our triple of
matrices reads
$$A=\left(\begin{array}{cc}
   1  & 1 \\
0     & 1
\end{array}\right), \quad B=\left(\begin{array}{cc}
    b_2 & 0 \\
    0 & 0
\end{array}\right), \quad C=\left(\begin{array}{cc}
    0 & c_1 \\
    c_1 & c_2
\end{array}\right).$$
For $\epsilon>0$ we consider the family of matrices
$$R_\epsilon=\left(\begin{array}{cc}
    \epsilon & 0 \\
    0 & \frac{1}{\epsilon^2}
\end{array}\right).$$
Acting with this family of matrices on the triple $(A,B,C)$ we obtain
\begin{eqnarray*}
R_\epsilon A R_\epsilon^{-1}&=& \left(\begin{array}{cc}
1 & \epsilon^3\\
0 & 1 \end{array}\right)\\
R_\epsilon B R_\epsilon^T&=&\left(\begin{array}{cc}
\epsilon^2 b_2 & 0\\
0 & 0\end{array}\right)\\
(R_\epsilon^T)^{-1} C R_\epsilon^{-1}&=&\left(\begin{array}{cc}
0 & \epsilon c_1\\
\epsilon c_1 & \epsilon^4 c_2\end{array}\right)
\end{eqnarray*}
and therefore
\begin{eqnarray*}
\lim_{\epsilon \to 0}R_\epsilon A R_\epsilon^{-1}&=& \left(\begin{array}{cc}
1 & 0\\
0 & 1 \end{array}\right)\\
\lim_{\epsilon \to 0}R_\epsilon B R_\epsilon^T&=&\left(\begin{array}{cc}
0 & 0\\
0 & 0\end{array}\right)\\
\lim_{\epsilon \to 0}(R_\epsilon^{T})^{-1} C R_\epsilon^{-1}&=&\left(\begin{array}{cc}
0 & 0\\
0 & 0\end{array}\right).
\end{eqnarray*}
We see that in this case $I_4$ lies in the closure of the
$\mathrm{GL}_2(\mathbb{R})$-orbit of $M_{A,B,C}$. The case where 
$c_1=0$ is analogous. One just needs to use the family of matrices
$$R_\epsilon=\left(\begin{array}{cc}
    \epsilon^2 & 0 \\
    0 & \frac{1}{\epsilon}
\end{array}\right).$$
It remains to discuss the case where $A$ is the identity matrix
$$A=I=\left(\begin{array}{cc}
    1 & 0 \\
    0 & 1
\end{array}\right).$$
The equation $A^2-BC=I$ implies in this case that
$$BC=0.$$
We first consider the case $B=0$. In this case the family
of matrices $\tfrac{1}{\epsilon}I$ acts on the triple $(A,B,C)$ by
$$\big(\tfrac{1}{\epsilon}I\big)_*(I,0,C)=(I,0,\epsilon^2 C)$$
with limit
$$\lim_{\epsilon \to 0}\big(\tfrac{1}{\epsilon}I\big)_*(I,0,C)=(I,0,0)$$
which shows that in this case $I_4$ lies in the closure of
the $\mathrm{GL}_2(\mathbb{R})$-orbit $M_{A,B,C}$. A similar argument
holds in the case where $C=0$ by using the family of matrices
$\epsilon I$ instead. It remains to discuss the case where neither
$B$ nor $C$ are the zero-matrix. Since $B$ is symmetric and transforms
as a bilinear form we can diagonalize $B$ so that we can assume without
loss of generality that
$$B=\left(\begin{array}{cc}
    b_1 & 0 \\
    0 & b_2
\end{array}\right)$$
with $b_1 \neq 0$. Since $C$ is symmetric as well we obtain from the
equation that $C$ has to be of the form
$$C=\left(\begin{array}{cc}
    0 & 0 \\
    0 & c_2
\end{array}\right).$$
Since $C$ is not the zero-matrix we must have $c_2 \neq 0$ which implies
in view of $BC=0$ that $b_2=0$ so that our triple becomes
$$A=\left(\begin{array}{cc}
    1 & 0 \\
    0 & 1
\end{array}\right), \quad B=\left(\begin{array}{cc}
    b_1 & 0 \\
    0 & 0
\end{array}\right), \quad C=\left(\begin{array}{cc}
    0 & 0 \\
    0 & c_2
\end{array}\right).$$
We consider the family of matrices
$$R_\epsilon=\left(\begin{array}{cc}
    \epsilon & 0 \\
    0 & \frac{1}{\epsilon}
\end{array}\right).$$
We act with this family of matrices on the triple $(A,B,C)$ to obtain
\begin{eqnarray*}
R_\epsilon A R_\epsilon^{-1}&=&A\\
R_\epsilon B R_\epsilon^T&=&\left(\begin{array}{cc}
\epsilon^2 b_1 & 0\\
0 & 0\end{array}\right)\\
(R_\epsilon^T)^{-1} C R_\epsilon^{-1}&=&\left(\begin{array}{cc}
0 & 0\\
0 & \epsilon^2 c_2\end{array}\right)
\end{eqnarray*}
and hence
\begin{eqnarray*}
\lim_{\epsilon \to 0}R_\epsilon A R_\epsilon^{-1}&=& \left(\begin{array}{cc}
1 & 0\\
0 & 1 \end{array}\right)\\
\lim_{\epsilon \to 0}R_\epsilon B R_\epsilon^T&=&\left(\begin{array}{cc}
0 & 0\\
0 & 0\end{array}\right)\\
\lim_{\epsilon \to 0}(R_\epsilon^{T})^{-1} C R_\epsilon^{-1}&=&\left(\begin{array}{cc}
0 & 0\\
0 & 0\end{array}\right).
\end{eqnarray*}

\begin{figure}
    \centering
    \includegraphics[width=0.99 \linewidth]{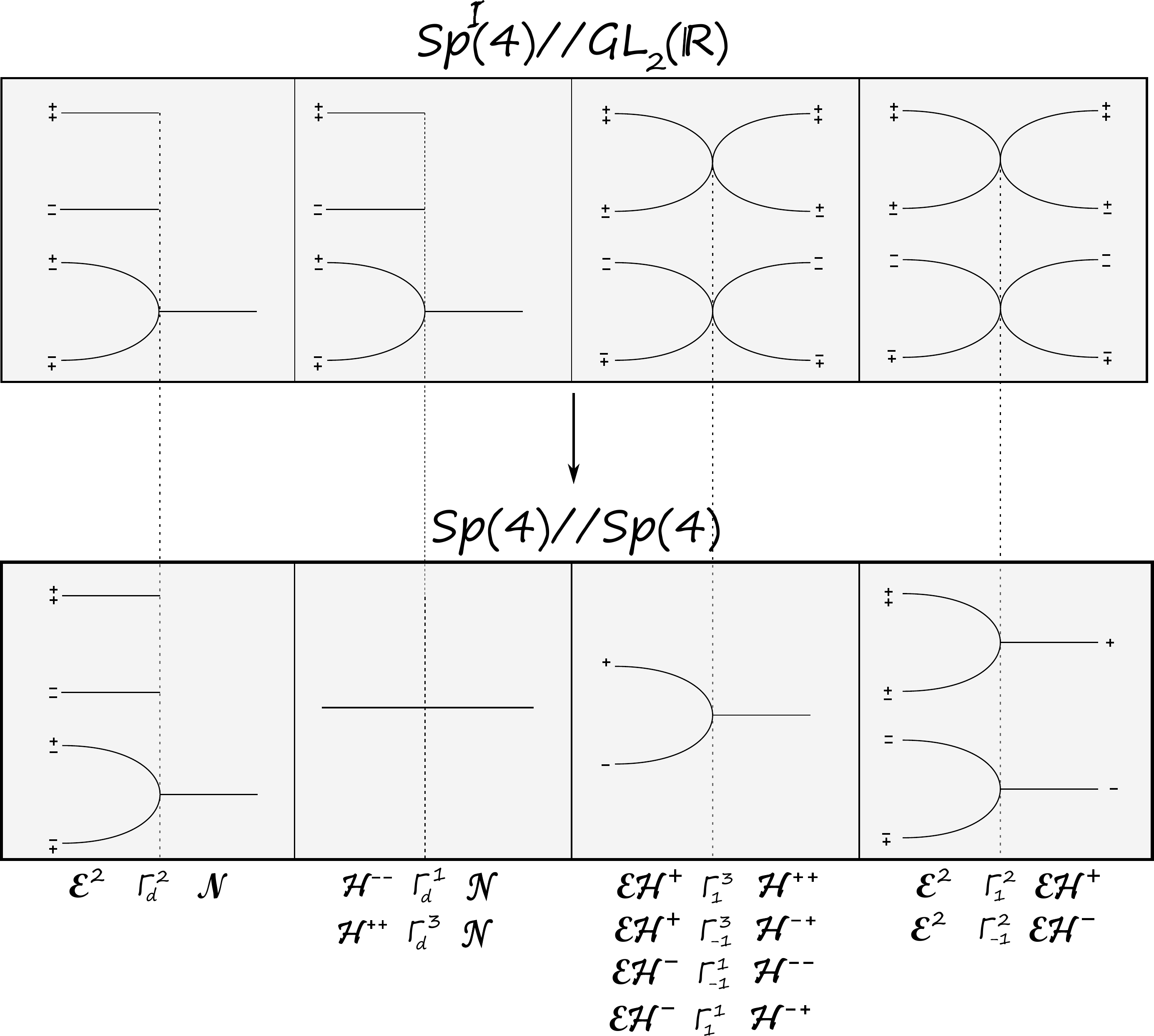}
    \caption{Different ($2$-dimensional) branches of $\mathrm{Sp}^\mathcal{I}(4)//\mathrm{GL}_2(\mathbb{R})$ and $\mathrm{Sp}(4)//\mathrm{Sp}(4)$. The signs on each branch correspond to $B$-positivity/negativity of the corresponding eigenvalue (a priori there are $4$ possibilities). The first vignette shows how they come together when crossing from $\mathcal{E}^2$ to $\mathcal{N}$ along $\Gamma_d^2$. On the second, when crossing from $\mathcal{H}^{--}$ to $\mathcal{N}$ along $\Gamma_d^1$; the picture is the same for $\mathcal{H}^{++}$ to $\mathcal{N}$ along $\Gamma_d^3$, and so on. All branches come together to a single point along each of the three singular points $(2,1),(0,-1), (-2,1)$. The natural map $\mathrm{Sp}^\mathcal{I}(4)//\mathrm{GL}_2(\mathbb{R})\rightarrow \mathrm{Sp}(4)//\mathrm{Sp}(4)$ in the GIT sequence collapses branches together, as shown in the picture. For example, above, $B$-positivity/negativity over the hyperbolic eigenspace of matrices of type $\mathcal{EH}^+$ is not invariant under symplectic conjugation, and hence the corresponding branches come together in $\mathrm{Sp}(4)//\mathrm{Sp}(4)$; we may still distinguish the elliptic eigenvalues via Krein theory, however.}
    \label{fig:bifurcations}
\end{figure}

This proves that as well in this last case $I_4$ lies in the closure
of the $\mathrm{GL}_2(\mathbb{R})$-orbit of $M_{A,B,C}$ and therefore
over $(2,1)$ the branched covers $\mathrm{Sp}^\mathcal{I}(4)//\mathrm{GL}_2(\mathbb{R})$ as well as
$\mathrm{Sp}(4)//\mathrm{Sp}(4)$ consists of a single point namely
the equivalence class of the matrix $I_4$.
\\ \\
The story over $(-2,1)$ is completely analogous. Over this point the
two branched covers just consist of the equivalence class of the matrix
$-I_4$.
\\ \\
We are left with a last point, namely $(0,-1)$. If a matrix $A$ lies over
this point it has $1$ and $-1$ as eigenvalues. In particular, it is diagonalizable and hence after taking advantage of the $\mathrm{GL}_2(\mathbb{R})$-action we can assume that $A$ has the form
$$A=\left(\begin{array}{cc}
    1 &  0\\
    0 & -1
\end{array}\right).$$
In particular, we have $A=A^T$ implying that the symmetric matrices
$B$ and $C$ commute with $A$. In particular, they have to keep invariant
the eigenspaces of $A$ and are therefore themselves diagonal matrices
$$B=\left(\begin{array}{cc}
    b_1 & 0 \\
    0 & b_2
\end{array}\right), \quad C=\left(\begin{array}{cc}
    c_1 & 0 \\
    0 & c_2
\end{array}\right).$$
The formula $A^2-BC=I$ implies that
$$b_1 c_1=0, \qquad b_2 c_2=0,$$
i.e., $b_1$ or $c_1$ has to vanish or $b_2$ or $c_2$ has to vanish. 
By going over to the GIT-quotient we can arrange that all of them vanish.
For example if $b_1$ does not vanish, then $c_1$ has to vanish and we
use the sequence of matrices
$$R_\epsilon=\left(\begin{array}{cc}
    \epsilon & 0  \\
    0 & 1
\end{array}\right)$$
to arrange that in the limit as $\epsilon$ goes to zero $b_1$  vanishes as well. Similarly, if $c_1$ does not vanish, then $b_1$ has to vanish and in this case we use the sequence of matrices
$$R_\epsilon=\left(\begin{array}{cc}
    \frac{1}{\epsilon} & 0  \\
    0 & 1
\end{array}\right)$$
and similarly for $b_2$ and $c_2$. Hence we can assume that $B=C=0$ and therefore we have the unique normal form
$$\left(\begin{array}{cccc}
    1 & 0 & 0 & 0  \\
    0 & -1&0 & 0\\
    0 & 0 & 1 & 0\\
    0 & 0 & 0&-1
\end{array}\right).$$
In particular, the two branched covers consists as well of a single point
over $(0,-1)$. This finishes the description of normal forms in all cases.

\section{Bifurcations and stability}

Given a family of symmetric spatial orbits, one considers the linearized flow along them, which induces a family of matrices in $\mathrm{Sp}^\mathcal{I}(4)//\mathrm{GL}_2(\mathbb{R})$. Note that a bifurcation of the family of orbits corresponds to a crossing of the eigenvalue $1$ of the family of matrices, or of $-1$ (period doubling), which geometrically means that the family crosses the walls $\Gamma_{\pm 1}$ when projected to $M_{2\times 2}(\mathbb{R})//GL_2(\mathbb{R})$. 

Recall that linear stability of an orbit means that the eigenvalues of the linearized matrix have strictly negative real part; in the symplectic case, this is equivalent to all eigenvalues having norm $1$, due to the symmetries of the spectrum. Moreover, strong linear stability means linear stability, even after small perturbations, i.e.\ it is a robust version of linear stability; we review these definitions in Appendix \ref{app:Krein}. 

In our setup, the (linearly) stable orbits are the ones whose matrices lie over the $\mathcal{E}^2$ component of the GIT quotient $M_{2\times 2}(\mathbb{R})//GL_2(\mathbb{R})$. The strongly (linearly) stable orbits correspond to those lying over the interior of $\mathcal{E}^2$, i.e.\ they cannot be perturbed to lie away from $\mathcal{E}^2$. However, there are also matrices which are strongly stable and lying over $\Gamma_d^2$, corresponding to the boundary of the $++$ and $--$ branches (see Figure \ref{fig:bifurcations}). The relationship between linear stability and the diagram of Figure \ref{fig:bifurcations} was already observed in \cite{Broucke}. Combining Figure \ref{fig:branchlocus} and Figure \ref{fig:bifurcations}, a moment's thought shows that, if we remove the bifurcation locus (the corresponding preimages of $\Gamma_{\pm 1}$), we have 8 connected components of its complement in $\mathrm{Sp}(4)//\mathrm{Sp}(4)$, and 19 of its complement in $Sp^\mathcal{I}(4)//GL_2(\mathbb{R})$.

\appendix

\section{The GIT quotient} \label{app:GIT}

In this appendix, we review the definition of the GIT quotient, and some nice general facts about the GIT quotient corresponding to the conjugation action of the general linear group on the space of matrices. 

Assume that $G$ is a Lie group which acts on a manifold $X$. 
The space of orbits $X/G$ is in general not a Hausdorff space. In order to remedy this
situation in some cases we consider the \emph{orbit closure relation} on $X$, namely
$$x \sim y \quad :\Longleftrightarrow \quad \overline{Gx}\cap \overline{Gy} \neq \emptyset,$$
i.e., the closures of the orbits of $x$ and of $y$ intersect. This relation is obviously 
reflexive and symmetric. If it is in addition transitive, it is an equivalence relation and
in this case we define the GIT quotient as
$$X//G:=X/\sim.$$
The following example plays an important role in our story.
\begin{prop}
The group $\mathrm{GL}_n(\mathbb{R})$ acts on the space of real $n\times n$-matrices
$\mathrm{M}_{n\times n}(\mathbb{R})$ by conjugation
$$R_* A=RAR^{-1}, \quad R \in \mathrm{GL}_n(\mathbb{R}),\,\,A \in \mathrm{M}_{n\times n}(\mathbb{R}).$$
For this action the orbit closure relation is transitive and the GIT quotient
$\mathrm{M}_{n \times n}(\mathbb{R})//\mathrm{GL}_n(\mathbb{R})$ is homeomorphic to
$\mathbb{R}^n$. If for $A \in \mathrm{M}_{n\times n}(\mathbb{R})$ the characteristic polynomial
is written as
$$p_A(t)=(-1)^n t^n+c_{n-1}t^{n-1}+\ldots +c_1 t+c_0,$$
then a homeomorphism is given by mapping the equivalence class of a matrix to the coefficients
of its characteristic polynomial
$$\mathrm{M}_{n \times n}(\mathbb{R})//\mathrm{GL}_n(\mathbb{R})\to \mathbb{R}^n, \quad
[A] \mapsto (c_{n-1},\ldots,c_1,c_0).$$
\end{prop}
\textbf{Proof: } Suppose that $A \in \mathrm{M}_{n\times n}(\mathbb{R})$. Then $A$ is conjugated
by a matrix $R \in \mathrm{GL}_n(\mathbb{R})$ to a matrix in real Jordan form. In case all
eigenvalues of $A$ are real, the real Jordan form does not differ from the complex Jordan form.
In case $A$ has as well nonreal eigenvalues its real Jordan form does not agree with its complex
Jordan form and needs some explanation. We first note that nonreal eigenvalues of $A$ 
appear in pairs $(\lambda,\overline{\lambda})$, since the matrix $A$ is real. In order to avoid
double counting we restrict our attention to nonreal eigenvalues in the upper halfplane
$\mathbb{H}=\{z=x+iy \in \mathbb{C}: y>0\}$.
\\ \\
If $\lambda \in \mathbb{R}$ and $m \in \mathbb{N}$ we define the Jordan block 
$J_{\lambda,m}$ as in the complex case
as the $m\times m$-matrix whose diagonal entries are all $\lambda$, and whose superdiagonal entries are
all $1$, while all other entries are zero like. For example,
$$J_{\lambda,3}=\left(\begin{array}{ccc}
\lambda & 1 & 0\\
0 & \lambda & 1\\
0 & 0 & \lambda
\end{array}\right).$$
If $\lambda=a+bi \in \mathbb{H}$ we first define the $2\times 2$-matrix
$$C_\lambda=\left(\begin{array}{cc}
a & -b\\
b& a
\end{array}\right).$$
Then different from the complex case we define for $m \in \mathbb{N}$ the Jordan block
$J_{\lambda,m}$ as the $2m\times 2m$-matrix consisting of $m\times m$ blocks of
$2 \times 2$-matrices whose diagonal entries are all $C_\lambda$, and whose superdiagonal entries are all
$I_2$, i.e., the $2 \times 2$-identity matrix, while all other entries are zero. For instance we have
$$J_{\lambda,3}=\left(\begin{array}{ccc}
C_\lambda & I_2 & 0\\
0 & C_\lambda & I_2\\
0 & 0 & C_\lambda
\end{array}\right).$$
A real Jordan matrix is then as usual a block matrix having Jordan blocks on the diagonal and
zeros elsewhere. 
\\ \\
Each Jordan block is similar to one where the superdiagonal is scaled by $\epsilon>0$. We illustrate
this paradigmatically for the Jordan block $J_{\lambda,3}$ for real $\lambda$, for which we have
$$\left(\begin{array}{ccc}
\epsilon^2 & 0 & 0\\
0 & \epsilon & 0\\
0 & 0 & 1
\end{array}\right)
\left(\begin{array}{ccc}
\lambda & 1 & 0\\
0 & \lambda & 1\\
0 & 0 & \lambda
\end{array}\right)
\left(\begin{array}{ccc}
\frac{1}{\epsilon^2} & 0 & 0\\
0 & \frac{1}{\epsilon} & 0\\
0 & 0 & 1
\end{array}\right)=
\left(\begin{array}{ccc}
\lambda & \epsilon & 0\\
0 & \lambda & \epsilon\\
0 & 0 & \lambda
\end{array}\right).$$
In particular, in the orbit closure of the matrix $A$ there lies a block diagonal matrix, 
namely 
$$D_A=\bigoplus_{\lambda \in \mathfrak{S}(A)\cap \mathbb{R}}\lambda^{\oplus \mathfrak{a}(\lambda)}
\oplus \bigoplus_{\lambda \in \mathfrak{S}(A) \cap \mathbb{H}}C_\lambda^{\oplus \mathfrak{a}(\lambda)},$$
where $\mathfrak{S}(A)$ denotes the spectrum of $A$, i.e., the set of all eigenvalues of $A$,
and $\mathfrak{a}(\lambda)$ denotes the algebraic multiplicity of an eigenvalue. Strictly speaking
we need to specify an order on the eigenvalues, in order to make $D_A$ well-defined as a matrix.
We choose the lexicographic order with real value as the first letter and imaginary value as the second one. Since however different ordering conventions lead to conjugated matrices the reader is
free to choose his own preferred convention which will not influence the following arguments. 
We note that $D_A$ is uniquely determined by the characteristic polynomial of $A$.
In particular, we see that if two matrices $A$ and $B$ have the same characteristic polynomial
$p_A=p_B$ we have
$$D_A \in \overline{\mathrm{GL}_n(\mathbb{R}) A} \cap \overline{\mathrm{GL}_n(\mathbb{R}) B}$$
implying that $A \sim B$. On the other hand, suppose that $A \sim B$. This means that there
exists a matrix
$$D \in \overline{\mathrm{GL}_n(\mathbb{R}) A} \cap \overline{\mathrm{GL}_n(\mathbb{R}) B}.$$
In particular, there exist a sequence $R_\nu \in \mathrm{GL}_n(\mathbb{R})$ such that
$$\lim_{\nu \to \infty} R_\nu A R_\nu^{-1}=D$$
as well as a sequence $S_\nu \in \mathrm{GL}_n(\mathbb{R})$ with
$$\lim_{\nu \to \infty} S_\nu A S_\nu^{-1}=D.$$
Since conjugated matrices have the same characteristic polynomial we have
$p_A=p_{R_\nu A R_\nu^{-1}}$ for every $\nu$ and therefore
$$p_D=\lim_{\nu \to \infty} p_{R_\nu A R_\nu^{-1}}=p_A$$
and similarly 
$$p_D=p_B$$
implying that 
$$p_A=p_B.$$
We have proved that
$$A \sim B \quad \Longleftrightarrow \quad p_A=p_B.$$
Since every polynomial with leading coefficient $(-1)^n$ arises as the characteristic polynomial
of a matrix $A \in \mathrm{M}_{n \times n}(\mathbb{R})$, the proposition follows. \hfill $\square$

\section{Krein theory and strong stability for Hamiltonian systems}\label{app:Krein}

In this appendix, we review some basic facts about Krein theory, its relationship with stability for orbits of Hamiltonian systems, and compare it to our notion of $B$-positivity in the case of symmetric orbits. We follow the exposition in Ekeland's book \cite{Eke90} (see also Abbondandolo's book \cite{Abo}).

Consider a linear symplectic ODE
$$
\dot x=M(t)x,
$$
where $M(t)=JA(t) \in \mathfrak{sp}(2n)$, with $A(t)$ symmetric, and $T$-periodic, i.e.\ $A(t+T)=A(t)$ for all $t$, and $J=\left(\begin{array}{cc}
    0 &  I\\
    -I & 0
\end{array}\right)$ is the standard complex multiplication. The solutions are given by $x(t)=R(t)x(0)$, where
$R(t)\in \mathrm{Sp}(2n)$ is symplectic and solves $\dot R(t)=M(t)R(t)$, $R(0)=I$.

\smallskip

\begin{fed}\textbf{(stability)}
The ODE $\dot x=JA(t)x$ is called \emph{stable} if all solutions remain bounded for all $t\in \mathbb{R}$. It is \emph{strongly} stable if there exists $\epsilon>0$ such that, if $B(t)$ is symmetric and satisfies $\Vert A(t)- B(t)\Vert <\epsilon$, then the ODE $\dot x=JB(t)x$ is stable. Similarly, a symplectic matrix $R$ is \emph{stable} if all its iterates $R^k$ remain bounded for $k\in \mathbb{Z}$, and it is \emph{strongly} stable if there exists $\epsilon>0$ such that all symplectic matrices  $S$ with $\Vert R-S \Vert<\epsilon$ are also stable.
\end{fed} 

Appealing to Floquet theory, one can show that the ODE $\dot x=JA(t)x$ is (strongly) stable if and only if $R(T)$ is (strongly) stable; see \cite[Section 2, Proposition 3]{Eke90}. Moreover, stability is equivalent to $R(T)$ being diagonalizable (i.e.\ all eigenvalues are semi-simple, meaning that their algebraic and geometric multiplicities agree), with its spectrum lying in the unit circle \cite[Section 1, Proposition 1]{Eke90}. Questions about the strong stability of Hamiltonian systems are therefore reduced to questions about the strong stability of symplectic matrices.

Now, recall that the spectrum of a symplectic matrix $R$ satisfies special symmetries. Concretely, its eigenvalues come in families of the form $\{\lambda,\overline{\lambda},\lambda^{-1}, \overline{\lambda}^{-1}\}.$ Therefore, if $\pm 1$ are eigenvalues, then they have even multiplicity. Moreover, if all its eigenvalues are simple, different from $\pm 1$, and lie in the unit circle, then they come in pairs $\{\lambda,\overline{\lambda}\}$. In this case, this implies that any other symplectic matrix close to $R$ will also have simple eigenvalues in the unit circle different from $\pm 1$ (otherwise an eigenvalue would have to bifurcate into two, which is not possible if eigenspaces are $1$-dimensional). Therefore in this case, $R$ is strongly stable. The case of eigenvalues with higher multiplicity is the subject of Krein theory, which we now review.

Consider the nondegenerate bilinear form $G(x,y)=\langle -iJx,y\rangle$ on $\mathbb{C}^{2n}$, associated to the Hermitian matrix $-iJ$. Every real symplectic matrix $R$ preserves $G$. Moreover, if $\lambda,\mu$ are eigenvalues of $R$ which satisfy $\lambda \overline{\mu}\neq 1$, then the corresponding eigenspaces are $G$-orthogonal, since
$$
G(x,y)=G(Rx,Ry)=\lambda \overline{\mu}G(x,y), 
$$
if $x,y$ are the corresponding eigenvectors. Moreover, if we consider the generalized eigenspaces 
$$
E_\lambda=\bigcup_{m\geq 1}\ker(R-\lambda I)^m,
$$
then it also holds that $E_\lambda,E_\mu$ are $G$-orthogonal if $\lambda \overline{\mu}\neq 1$ \cite[Section 2, Proposition 5]{Eke90}. This, in particular, implies that if $\vert\lambda\vert \neq 1$, then $E_\lambda$ is $G$-isotropic, i.e.\ $G\vert_{E_\lambda}=0$. If $\sigma(R)$ denotes the spectrum of $R$, we have a $G$-orthogonal decomposition
$$
\mathbb{C}^{2n}=\bigoplus_{\substack{\lambda \in \sigma(R)\\\vert \lambda \vert\geq 1}}F_\lambda,
$$
where $F_\lambda=E_\lambda$ if $\vert \lambda\vert=1$, and $F_\lambda=E_\lambda \oplus E_{\overline{\lambda}^{-1}}$ if $\vert\lambda\vert >1$. Since $G$ is non-degenerate, and the above splitting is $G$-orthogonal, the restriction $G_\lambda=G\vert_{F_\lambda}$ is also non-degenerate. Recall that the signature of a non-degenerate bilinear form $G$ is the pair $(p,q)$, where $p$ is the dimension of a maximal subspace where $G$ is positive definite, and $q$ is the dimension of a maximal subspace where $G$ is negative definite. Note that if $\vert \lambda \vert \neq 1$, with algebraic multiplicity $d$, then the $2d$-dimensional space $F_\lambda$ has $E_\lambda$ as a $d$-dimensional isotropic subspace, and hence the signature of $G_\lambda$ is $(d,d)$. On the other hand, if $\vert\lambda\vert=1$, then the non-degenerate form $G_\lambda$ can have any signature. This justifies the following:
\begin{fed}\textbf{(Krein-positivity/negativity)}
If $\lambda$ is an eigenvalue of the symplectic matrix $R$ with $\vert \lambda \vert=1$, then the signature $(p,q)$ of $G_\lambda$ is called the \emph{Krein-type} or \emph{Krein signature} of $\lambda$. If $q=0$, i.e.\ $G_\lambda$ is positive definite, $\lambda$ is said to be \emph{Krein-positive.} If $p=0$, i.e.\ $G_\lambda$ is negative definite, $\lambda$ is said to be \emph{Krein-negative.} If $\lambda$ is either Krein-negative or Krein-positive, we say that it is \emph{Krein-definite}. Otherwise, we say that it is \emph{Krein-indefinite}.
\end{fed}

If $\lambda$ is of Krein-type $(p,q)$, then $\overline{\lambda}$ is of Krein-type $(q,p)$ \cite[Section 2, Lemma 9]{Eke90}. If $\lambda$ satisfies $\vert \lambda \vert=1$ and it is not semi-simple, then it is easy to show that it is Krein-indefinite \cite[Section 2, Proposition 7]{Eke90}. Moreover, $\pm 1$ are always Krein-indefinite if they are eigenvalues, as they have real eigenvectors $x$, which are therefore $G$-isotropic, i.e.\ $G(x,x)=0$. The following, originally proved by Krein in \cite{Kre1,Kre2,Kre3,Kre4} and independently rediscovered by Moser in \cite{Moser}, gives a characterization of strong stability in terms of Krein theory:
\begin{thm}\label{Kreinthm}
 $R$ is strongly stable if and only if it is stable and all its eigenvalues are Krein-definite. 
\end{thm}
See \cite[Section 2, Theorem 3]{Eke90} for a proof. Note that this generalizes the case where all eigenvalues are simple, different from  $\pm 1$ and in the unit circle, as discussed above.

\medskip

We now prove Lemma \ref{lemma:BposKrein}. 

\begin{proof}[Proof of Lemma \ref{lemma:BposKrein}] Since the notion of Krein-type
is invariant under symplectic conjugation and only involves the
eigenspaces it suffices to show that for the matrix
$$M=\left(\begin{array}{cc}
  \cos \theta   &  -\sin \theta\\
   \sin \theta  & \cos \theta
\end{array}\right)$$
the eigenvalue $e^{i\theta}$ is Krein-negative; the positive case is analogous. An eigenvector
is given by
$$v=\left(\begin{array}{c}
     1  \\
     -i
\end{array}\right).$$
We have
$$G(v,v)=-2$$
and this shows that $e^{i\theta}$ is Krein-negative, concluding the proof.\end{proof}

\end{document}